\documentclass[a4paper,10pt]{article}

\usepackage[utf8]{inputenc}   
\usepackage[T1]{fontenc}      
\usepackage{geometry}         
\usepackage[francais,english]{babel}  

\title{The integrability of the $2$-Toda lattice on simple Lie algebra}           
\author{Khaoula Ben Abdeljelil}
\date{}                       

\newtheorem{theorem}{Theorem}
\newtheorem{proposition}[theorem]{Proposition}
\newtheorem{lemma}[theorem]{Lemma}
\newtheorem{corollary}[theorem]{Corollary}

\newtheorem{conjecture}[theorem]{Conjecture}
\newtheorem{definition}[theorem]{Definition}

\newtheorem{example}[theorem]{Example}

\newtheorem{remark}[theorem]{Remark}

\newcommand{\qed}{\nobreak \ifvmode \relax \else
      \ifdim\lastskip<1.5em \hskip-\lastskip
      \hskip1.5em plus0em minus0.5em \fi \nobreak
      \vrule height0.75em width0.5em depth0.25em\fi}

\makeindex 

\usepackage{pb-diagram}  
\usepackage{makeidx}     
\usepackage{graphicx}    
\usepackage{multicol}    
\usepackage{cite}        
\usepackage{amsfonts}
\usepackage{amssymb}
\usepackage{amscd}
\usepackage{epsfig}

\def\C{\mathbf{C}}

\def\Z{\mathbf{Z}}
\def\R{\mathbf{R}}
\def\N{\mathbf{N}}

\def\a{\alpha}
\def\b{\beta}

\def\l{\lambda}

\def\rmi{\uppercase\expandafter{\romannumeral1}}
\def\rmii{\uppercase\expandafter{\romannumeral2}}
\def\rmiii{\uppercase\expandafter{\romannumeral3}}
\def\rmv{\uppercase\expandafter{\romannumeral5}}
\def\rmvi{\uppercase\expandafter{\romannumeral6}}
\def\rmvii{\uppercase\expandafter{\romannumeral7}}
\def\rmviii{\uppercase\expandafter{\romannumeral8}}

\def\RR{\mathcal{R}}              

\def\FF{\mathcal{F}}          
\def\TT{\mathcal{T}}

\def\I{\mathcal{I}}

\def\X{\mathcal{X}}              

\def\g{\mathfrak{g}}          

\def\h{\mathfrak{h}}          
\def\Liesl#1{\mathfrak{sl}_{#1}}     
\def\gl#1{\mathfrak{gl}_{#1}}        

\def\can#1{\left\langle#1\right\rangle}

\def\PB{\left\{\cdot\,,\cdot\right\}}

\def\pb#1{\left\{#1\right\}}

\def\lb#1{\[#1\]}
\def\LB{[\cdot\,,\cdot]}

\def\inn#1#2{\left\langle#1\,\vert\,#2\right\rangle}

\def\INN{\langle\cdot\,\vert\cdot\,\rangle}

\def\({\left(}
\def\){\right)}
\def\[{\left[}
\def\]{\right]}

\def\ad{\mathop{\rm ad}\nolimits}
\def\Ad{\mathop{\rm Ad}\nolimits}

\def\det{\mathop{\rm det}\nolimits}

\def\dim{\mathop{\rm dim}\nolimits}
\def\diff{\mathsf{d}}

\newcommand{\Rk}{\mathop{\rm Rk}\nolimits}

\def\card{\mathop{\rm card}\nolimits}

\newenvironment{eqn*}[1][1.5]
  {$$\renewcommand{\arraystretch}{#1}
      \begin{array}{rcl}}
      {\end{array}$$}

\def\comment#1{}  

\def\p{\partial}
\def\pp#1#2{\frac{\p #1}{\p #2}}

\begin{document}

\author{Khaoula Ben Abdeljelil}
\title{The integrability of the periodic Full Kostant-Toda  on  a simple Lie algebra}
  \date{}                     
\maketitle

\begin{abstract}
 We define the periodic Full Kostant-Toda  on every simple Lie algebra,  and show its Liouville integrability. More 
 precisely we  show that this  lattice is given by a Hamiltonian vector field, associated to a Poisson bracket which results from an $\RR$-matrix. We construct a large family of constant of motion which we  use to prove the Liouville integrability of the system with the help of several results on simple Lie algebras, $\RR$-matrix, invariant functions and root systems.    
\end{abstract}
 \tableofcontents
\section{Introduction}\label{sect0}
The \emph{non-periodic (resp. periodic) Toda lattice on $\Liesl{n} (\C)$} is the system of differential equations 
given by a following Lax equation:
 \begin{equation}\label{toda-lax-equation}
{\dot{L}}=\lb{L,L_-}, \quad  \quad (\mbox{resp. }  {\dot{L}(\l)}=\lb{L(\l),L(\l)_-} ),
\end{equation}%
where $L$ and $L_-$ are the traceless matrices of the form given below.
For the non-periodic case, we impose:
\begin{equation}\label{todanp}
L=\left(
\renewcommand{\arraystretch}{0.9}
\begin{array}{cccccc}
b_1&1&0&\cdots&\cdots &0\\
a_1&b_2&1&\ddots&&\vdots \\
0&\ddots&\ddots&\ddots&\ddots&\vdots\\
\vdots&\ddots&\ddots&\ddots&\ddots&0\\
 \vdots&  &\ddots&a_{n-2}&b_{n-1}&1\\
0&\cdots&\cdots&0 &a_{n-1}&b_n
\end{array}
\right)
,
L_-=\left(
\renewcommand{\arraystretch}{0.5}
\begin{array}{cccccc}
0     &      &     0&\cdots&\cdots &0\\
a_1   &0     &      &\ddots& &\vdots\\
0     &\ddots&\ddots&\ddots  &\ddots&\vdots\\
\vdots&\ddots&\ddots&\ddots  &\ddots&0\\
\vdots&      &\ddots&a_{n-2}&0&0\\
     0&\cdots&\cdots&0 &a_{n-1}&0
\end{array}
\right).
\end{equation}%
In the periodic case, we choose a formal parameter $\l$ and we impose:
\begin{equation}\label{a32}
L(\l)=\left(
\renewcommand{\arraystretch}{0.5}
\begin{array}{cccccc}
b_1&1&0&\cdots&0 &a_n{\lambda}^{-1}\\
a_1&b_2&1&\ddots&&0 \\
0&\ddots&\ddots&\ddots&\ddots&\vdots\\
\vdots&\ddots&\ddots&\ddots&\ddots&0\\
 0&  &\ddots&a_{n-2}&b_{n-1}&1\\
\lambda&0&\cdots&0 &a_{n-1}&b_n
\end{array}
\right), \quad
L(\l)_-=
\left(
\renewcommand{\arraystretch}{0.5}
\begin{array}{ccccc}
0   & \cdots&\cdots&0&    a_n{\lambda}^{-1}\\
a_1&0     &   &  &0\\
0&\ddots&\ddots&&\vdots\\
\vdots&\ddots&a_{n-2}&\ddots&\vdots\\
0&\cdots&0&a_{n-1}&0
\end{array}
\right).
\end{equation}%
These systems of differential equations are classical examples of what is called Liouville
 integrable systems \cite[Definition 4.13]{Liv},
which form a class of equations known to be integrable by quadrature (i.e.,
 whose solutions can be expressed from their initial values with the help of elementary operations,
 integration, and inversion of diffeomorphism, see \cite[Section 4.2]{Liv}  for a more precise description).
 For our present purpose, we have to introduce Liouville integrability not only for symplectic manifolds,
 but in the enlarged context of Poisson manifolds (see again \cite{Liv} for the notion 
of Poisson manifold, and related notions, like rank, Casimir functions and involutive families):
\begin{definition}\label{def:liouville}
Let $(M,\PB)$ be a Poisson manifold of rank $2r$. A family $\FF=(F_1,\dots,F_s)$ of functions on $M$
is said to be \emph{Liouville integrable} if
\begin{enumerate}
\item[(1)] For all $i,j=1,\dots, s$, the functions $F_i,F_j$ commute, i.e., $\pb{F_i,F_j}=0$.
\item[(2)] The functions  $(F_1,\dots,F_s)$ form an independent family on $M$.
\item[(3)] $s=\dim M-r$, i.e., $\card \FF=\dim M-\frac{1}{2}\Rk(M,\PB)$.
\end{enumerate}
The triple $(M,\PB,\FF)$ is then said  to be a \emph{Liouville  integrable system} of \emph{rank} $ 2r$.
\end{definition}
\noindent
By a slight abuse of vocabulary, a differential equation is said to be
Liouville integrable when one can find a Liouville integrable
system such that one of the Hamiltonian vector fields describes the equation.

The non-periodic and periodic Toda lattices admit a natural extension\footnote{Here, we do not wish to give a precise meaning to the word "extension", that we simply to use to speak of a differential equation of the same shape on a bigger phase space.} 
and several of them have been proved to be Liouville integrable.

To start with, Deift, Li, Nanda, Tomei  \cite{DLNT} have proved the Liouville integrability of
 the \emph{(non-periodic) Full Kostant-Toda lattice}, that they define to be the system of differential equations given by:
\begin{equation}\label{ft}
\dot L=[L,L_-],
\end{equation}%
where $L$ is a symmetric matrix of $\gl{n}(\C)$ and $L_-$  it is the  skew-symmetric part 
of $L$ with respect to the decomposition of matrices as upper-triangular matrices and skew-triangular matrices. Up to a Poisson morphism, this system is shown by Ercolani, Flaschka and Singer \cite{EFS} to be given by an equation of the form (\ref{ft}), where $L$ is of the form:
\begin{equation}\label{Full}
L=\left(
\renewcommand{\arraystretch}{0.5}
\begin{array}{cccc}
a_{11}&1     &         &0     \\
a_{21}&a_{22}&\ddots   &      \\
\vdots&      &\ddots   &1     \\
a_{n1}&\cdots&a_{n,n-1}&a_{nn}
\end{array}
\right)
\in\gl{n}(\C)
\end{equation}%
and $L_-$ is the strictly lower triangular part of $L$ with respect to the decomposition of matrices as upper-triangular matrices and strictly lower-symmetric matrices.\\
Exactly as the non-periodic Full Kostant-Toda lattice is an extension of the non-periodic Toda lattice, there is a natural extension
of the periodic Toda lattice, namely the system of differential equations is given by:
\begin{equation}\label{mpf}
\dot { L}(\lambda)=[L(\lambda),{ L(\lambda)}_-],
\end{equation}%
where $\l$ is a formal parameter and $L(\l)$ is imposed to be of the form:
\begin{equation}\label{Matft}
 L(\lambda)=\left(
 \renewcommand{\arraystretch}{0.5}
\begin{array}{ccccc}
a_{11}        &1+b_{12}\lambda ^{-1}                &b_{13}\lambda ^{-1}     &\cdots        &b_{1n}\lambda ^{-1}     \\
\vdots                 & \ddots ~~~~~~~~              &\ddots ~~~~~~               &\ddots        &\vdots      \\
\vdots                 &                            &\ddots ~~~~~~~~           &\ddots ~~~~~~    &b_{n-2,n}\lambda ^{-1}\\
a_{n-1,1}              & \cdots                     &\cdots                  &\ddots ~~~~~~~~ &1+b_{n-1,n} \lambda ^{-1}\\
a_{n1}+\lambda         & a_{n2}                     &\cdots                  &\cdots        &a_{nn}~~~~
\end{array}
\right)
\end{equation}%
and
 \begin{eqn*}
{ L(\lambda)}_-=\left(
\renewcommand{\arraystretch}{0.5}
\begin{array}{cccc}
      0        &b_{12}\lambda ^{-1} &\cdots     &b_{1n}\lambda ^{-1}     \\
a_{21}         & \ddots             &\ddots     &\vdots      \\
\vdots         & \ddots             &\ddots     &b_{n-1,n}\lambda ^{-1}\\
a_{n1}         &\cdots              &a_{n,n-1}  &0
\end{array}
\right).
\end{eqn*}%
We call this system of differential equations the \emph{periodic Full Kostant-Toda lattice on $\Liesl{n}(\C)$}.
It deserves to be noticed that this system (more precisely its symmetric equivalent)
 appears as the extreme case of a sequence of systems studied by van Moerbeke and Mumford \cite{MMD}.

For all the systems previously introduced on $ \Liesl{n}(\C) $, there is natural
 manner to replace $\Liesl{n}(\C)$ by an arbitrary simple Lie algebra $\g$. 
Liouville integrability has been proved for an arbitrary simple Lie algebra 
in the cases of the periodic and non-periodic Toda lattices see \cite{kostant}, and in 
the case of non-periodic Full Kostant-Toda lattice by Gekhtman and Shapiro \cite{Gekhtman}.
The purpose of the present article is to show the Liouville integrability of the periodic
 Full Kostant-Toda lattice for every simple Lie
algebra.
This article is organized as follows. To start with, we define the periodic Full Kostant-Toda lattice 
and its phase space for every simple Lie algebra $\g$ in Section \ref{sect1}. More precisely, 
we construct this space as a finite dimensional  affine subspace of the loop algebra $\g[\l,\l^{-1}]$.  
This phase space is endowed with a Poisson structure in Section \ref{sect2}.
 A celebrated theorem, called the AKS theorem (see \cite[Theorem 4.37]{Liv}), implies
 that all the coefficients in $\l$ of the  $\ad$-invariant
 functions on $\g[\l,\l^{-1}]$ commute, therefore this family is a good candidate to prove Liouville integrability.
In Section \ref{sect3}, by restricting this family to the phase space of the periodic Full Kostant-Toda lattice, we state
 the main theorem: the integrability of the periodic Full Kostant-Toda lattice on $\g$, 
the proof of which will
be separated in several steps. The independence of the family of functions that we consider will
 be proved in Proposition~\ref{funinde}, with a help of a sophisticated result about regular $\Liesl{2} (\C)$-triplets 
and $\ad$-invariant functions established by Ra\"is \cite{Rais}. But the most difficult point is the computation
 of the rank of the Poisson structure on $\TT_{\l}$. This computation will be done with the help of Maple
 for the exceptional simple Lie algebras and the treatment of the four series of regular simple Lie algebra
 is completed with the help of a detailled  investigation of the root system of those. 
In section \ref{sect4}, we
 finish this study by presenting a conjectured generalization.
\section{Definition of the periodic Full Kostant-Toda lattice on a  simple Lie algebra}\label{sect1}
 In this section, we define the $2$-Toda lattice on every  simple Lie algebra.\\
 Let $\g$ be  a simple Lie algebra of
rank $\ell$, with Killing form  $\INN$. We choose  $\h$ a
Cartan subalgebra with root system $\Phi$, and
$\Pi=(\a_1,\dots,\a_{\ell})$  a system of simple roots with
respect to $\h$. For every $\a$ in $\Phi\backslash \{-\Pi,\Pi\}$, we denote by  $e_{\a}$  a
non-zero eigenvector associated to eigenvalue $\a$ and, for every
$1\leqslant i\leqslant \ell$,  we denote  by $e_i$ and $e_{-i}$ a non-zero  eigenvector
associated respectively  to $\a_i$ and $-\a_i$.  The Lie algebra $\g=\bigoplus_{k\in\Z}\g_k$  is  endowed with the natural
 grading (i.e., for every $k,l\in\Z$, $[\g_k,\g_l]\subset\g_{k+l}$)
  defined by $\g_0:=\h$ and,  for every  $k\in\Z$,  $\g_k:=<e_{\a}\mid   \a\in\Phi,
  |\a|=k>$, where $|\a|$ is the  length of the root $\a$, i.e., $|\a|:=\sum_{i=1}^{\ell}a_i$ for  $\a=\sum_{i=1}^{\ell}a_i\a_i$ and we denote by $\b$ the longest root of $\g$.
Recall that: $\inn{\g_k}{\g_l}=0$ if $k+l\neq 0$.
 We introduce the following notation
\begin{eqn*}
\g_{<k}:=\bigoplus_{i<k}^{} \g _i,&\qquad \g_{\leqslant k}:=\bigoplus_{i\leqslant k}^{}\g _i,\\
\g_{>k}:=\bigoplus_{i>k}^{} \g _i,&\qquad \g_{\geqslant k}:=\bigoplus_{i\geqslant k}^{} \g_i.
\end{eqn*}%
 The next definition gives back  the  definition given in Section \ref{sect0} of the periodic
  Full Kostant-Toda lattice on $\Liesl{n}(\C)$
 when specialized to the case of  $\g=\Liesl{n}(\C)$ and $\h$ is a Lie  subalgebra formed by the diagonal 
matrices of $\Liesl{n}(\C)$.
 \begin{definition}
 The \emph{periodic Full Kostant-Toda lattice, associated to a simple Lie algebra~$\g$},
 is the  system of differential equations given by  the following Lax equation:
\begin{equation}\label{mvpf}
\dot { L}(\lambda)=[L(\lambda),{ L(\lambda)}_-],
\end{equation}%
where $L(\l)=\lambda e_{-\beta}+\sum_{i=1}^{\ell}(a_ih_i+e_i)+
\sum_{\a\in\Phi _+}^{}(a_{-\a}e_{-\a}+\lambda^{-1}b_{\a}e_{\a})$ is  an element of the following
phase space $\TT_{\l}$ of the  periodic
Full Kostant-Toda  lattice
\begin{equation}
\TT_{\lambda}:=\lambda^{-1}\g_{> 0}+(\g_{\leqslant 0}
+\sum_{i=1}^{\ell}e_i)+\lambda e_{-\beta}
\end{equation}%
 and
$L(\lambda)_-=\sum_{\a\in\Phi_+}(a_{-\a}e_{-\a}+\lambda^{-1}b_{\a}e_{\a})$.
\end{definition}
\section{Poisson structure on the phase space of the periodic Full Kostant-Toda lattice}
\label{sect2}
In the present section, we show that the periodic Full Kostant-Toda lattice is a Hamiltonian system, with respect to a Poisson
structure on $\TT_{\l}$, naturally obtained  as a substructure of a linear Poisson on the loop algebra $\g\otimes\C[\l,\l^{-1}]$, associated to an $R$-matrix.
\subsection{Poisson structure on the loop algebra $\g\otimes\C[\l,\l^{-1}]$}
Let   $\tilde\g$  the  loop algebra,
namely the tensor product
$\tilde\g=\g\otimes\C[\lambda,{\lambda}^{-1}]$, whose elements are sums
\begin{eqn*}
x(\lambda)=\sum_{i\in\Z}^{}x_i\lambda^i,
\end{eqn*}%
where finitely many $(x_i)_{i\in\Z}$ are non zero.
 We first endow  $\tilde \g$
  with the unique  bilinear  bracket   $\C[\lambda,\lambda^{-1}]$, which  extends the Lie bracket of  $(\g,\LB)$.\\
  We construct  a Poisson structure on
 the algebra of functions defined on the phase space of the periodic Full Kostant-Toda lattice.\\
We introduce a grading on
 $\tilde\g$ by  defining
the degree of  ${\lambda}^ke_{\alpha}$, ($\a$ being  a root of  $\g$ and
$k\in\Z$)  to be    $|\alpha|+(|\beta|+1)k$, where we recall  that  $\b$  is the longest positive root of $\g$.\\
We denote     by  ${\tilde{\g}}_i$ the Lie  subspace of weight
$i$, which  defined by:
\begin{eqn*}
{\tilde{\g}}_i:=\langle{{\lambda}^k e_{\alpha} \textrm{ such that } |\alpha|+(|\beta|+1)k=i,
\textrm{ for every  }  \alpha\in\Phi, k\in\Z}\rangle .
\end{eqn*}%
\begin{lemma}\label{juhb}~\
\emph{(1)}   For  $i=0$,  $\tilde {\g}_0=\h$, for every    $i=-|\beta|,\dots, -1 $, $
{\tilde\g}_i=\g_i\oplus \l^{-1}\g_{i+|\beta|+1}$
and for every  $i=1,\dots,|\beta|$, $
  {\tilde\g}_i=\g_i\oplus \l\g_{i-|\beta|-1}$.\\
\emph{(2)}  $\tilde\g=\bigoplus_{k\in\Z}{\tilde\g}_k$ is a graded Lie algebra and 
$\tilde{\g}$ admits the following vector space  decomposition:
\begin{equation}\label{spl-loop-algebra}
\tilde{\g}={\tilde{\g}}_+\oplus {\tilde{\g}}_-,
\end{equation}%
where
\begin{eqn*}
 {\tilde{\g}}_+:=\bigoplus_{i\geq{0}}{\tilde{\g}}_i\qquad\textrm{  and }\qquad {\tilde{\g}}_-:=\bigoplus_{i<0}{\tilde{\g}}_i
\end{eqn*}%
 are  Lie subalgebras of   $\tilde \g$.
\end{lemma}
Let  ${\tilde \g}^*$ be  the space of all linear   forms on $\tilde\g$  which are identically   zero on all $({\tilde\g}_i)_{i\in\Z}$
except  finitely many of them. We notice  that   the space $\tilde{\g}^*$
 have the following   decomposition:
\begin{eqn*}
 {\tilde{\g}}^*:= \bigoplus_{i\in\Z}{\tilde{\g}}^*_i,
\end{eqn*}%
where
\begin{eqn*}
{\tilde{\g}}^*_i:=\{\xi\in{\tilde{\g}}^*\mid \xi\textrm{ is zero on }\tilde{\g}_j,\textrm{ for every   }j\neq i\}.
\end{eqn*}%
 Let $\INN_{\l}$ be the following   non-degenerate, $\ad$-invariant, symmetric form:
\begin{equation}\label{isal}
\renewcommand{\arraystretch}{0.9}
\begin{array}{cccc}
\INN_{\lambda} :&\tilde\g\times\tilde\g&\to&\C\\
    &  (X(\lambda),Y(\lambda))&\mapsto&\sum_{k\in\Z}\inn{X_k}{Y_{-k}}.
\end{array}
\end{equation}%
  The bilinear    form (\ref{isal}) gives an identification between
$\tilde{\g}_i^*$ and    $\tilde{\g}_{-i}$, hence
between    ${\tilde \g}$ and   ${\tilde \g}^*$. Moreover,
 the orthogonal of  ${\tilde \g}_i$, for every  $i\in\Z$, is  ${\tilde
    \g}^{\perp}_i:=\bigoplus_{j\neq i}^{}{\tilde \g}_{-j}$.\\
Let
 $\FF(\tilde \g)$ be the symmetric algebra  generated by the elements of  ${\tilde \g}^*$
 (is a subalgebra of  the algebra of polynomial functions
   on  $\tilde \g$ and by  construction is such that the    gradient of a
 function in  a point of  $\tilde \g$  is   in   $\tilde \g$).  Then
   $\tilde\g$ is equipped of  the Poisson  structure\footnote{because ${\tilde \g}^*$ is equipped of the Poisson
${\tilde R}$-bracket and ${\tilde \g}^* \sim\tilde\g$.}, given
for every   $F,G\in\FF(\tilde\g)$ and every   $x(\lambda)\in\tilde\g$, by:
\begin{equation}\label{fdcv}
{\pb{F,G}}_{\tilde
  R}(x(\lambda))={\inn{x(\lambda)}{[\nabla_{x(\lambda)}F,\nabla_{x(\lambda)}G]_{\tilde
    R}}}_{\l},
\end{equation}%
 where  $\tilde R$ is an $\tilde R$-matrix of $\tilde\g$, defined by:
\begin{equation}\label{kl2ml}
\tilde R:={\tilde P}_+-{\tilde P}_-,
\end{equation}%
 and  ${\tilde P}_{\pm}$ is  the projection of  $\tilde \g$ on
${\tilde\g}_{\pm}$. For every  element $x(\l)$, we denote $x(\l)_{\pm}:={\tilde P}_{\pm}(x(\l))$.
  In formula (\ref{fdcv}),   $\nabla_{x(\lambda)}F$ stands for the  gradient of  $F$ at the point  $x(\lambda)$ computed
with respect to  ${\INN}_{\l}$.
\subsection{The Poisson ${\tilde R}$-bracket on $\FF(\TT_{\l})$}
The next proposition should be interpreted as meaning that $\TT_\l$ is a Poisson submanifold of
 $(\tilde{\g},\PB_{\tilde R})$, but the fact that $\tilde{g}$ is infinite dimensional prevents us to state it in that manner. What makes sense however is to show that there exists a unique Poisson bracket on
the algebra $\FF(\TT_\l)$ such that the restriction map $\FF(\tilde{\g})$ is a Poisson morphism. Indeed, since this restriction map is surjective, to prove the existence of this Poisson structure, it suffices to prove that the ideal $\I=\langle F\in\FF(\tilde\g) \mid F\equiv 0 \textrm{ on }\TT_{\l}\rangle$ is a Poisson ideal of the Poisson algebra $(\FF(\tilde{\g}), \PB_{\tilde{R}})$.
\begin{proposition}\label{fgyu}
 The phase space of the  periodic  Full Kostant-Toda
 $\TT_{\l}$  inherits an unique Poisson
  structure $(\FF(\tilde \g),{\PB}_{\tilde R})$
  such that the restriction map $\FF(\tilde{\g}) \to \FF( \TT_\l)$ is a Poisson morphism.
\end{proposition}
\begin{proof}
As stated before the proposition, we are left with the task of verifying that the ideal $\I$
 is a Poisson ideal with respect to the Poisson bracket~${\PB}_{\tilde R}$. According to Lemma~\ref{juhb}, the  affine subspace  $\TT_{\l}$ of ${\tilde\g}$ can be described as follows:
\begin{equation}\label{rel_ps}
\TT_{\lambda}:=\bigoplus_{-|\beta|\leqslant i\leqslant 0}{\tilde{\g}}_i +f,
\end{equation}%
where  $f:=\sum_{i=1}^{\ell}e_i+\lambda e_{-\beta}\in{\tilde \g}_1$.
The gradient at a point  $L(\l)\in\TT_{\l}$ of an arbitrary function $F\in\I$ satisfy the following relation:
\begin{equation}
\nabla_{L(\l)}F\in \bigoplus_{-|\b|\leqslant i\leqslant 0}^{}{\tilde
    \g}_i^{\perp}={\tilde \g}_{<0}\oplus{\tilde\g}_{\geq |\b|+1},
\end{equation}%
so that there exists  $x(\l) \in{\tilde \g}_{<0}$ and
$y(\l)\in{\tilde\g}_{\geq |\b|+1}$, such that  $\nabla_{L(\l)}F=x(\l)+y(\l)$.
For an arbitrary function  $G\in\FF(\tilde \g)$,
\begin{eqnarray*}
{\pb{F,G}}_{\tilde
  R}(L(\l))&=&\inn{L(\l)}{[(\nabla_{L(\l)}F)_+,(\nabla_{L(\l)}G)_+]-[(\nabla_{L(\l)}F)_-,(\nabla_{L(\l)}G)_-]}\\
          &=&\inn{L(\l)}{[y(\l),(\nabla_{L(\l)}G)_+]-[x(\l),(\nabla_{L(\l)}G)_-]}\\
          &=& 0,
\end{eqnarray*}%
where, in the last line, we have used the fact that  $L(\l)\in\bigoplus_{-|\b|\leqslant i\leqslant 1}^{}{\tilde
    \g}_i $ is  orthogonal to both $[y(\l),(\nabla_{L(\l)}G)_+]$ (which belongs to ${\tilde
  \g}_{\geq |\b|+1} $) and $[x(\l),(\nabla_{L(\l)}G)_-]$ (which belongs to ${\tilde
  \g}_{<-1}$). The ideal $\I$ is then a Poisson ideal, which endows to
 $(\FF(\tilde\g) /\I,{\PB}_{\tilde R})$ with a  Poisson $\tilde R$-bracket. Since the algebra  $\FF(\tilde\g) /\I$
  is  canonically isomorphic  to  $\FF(\TT_{\l})$, this Poisson $\tilde R$-bracket is an algebraic Poisson structure on
 $\TT_{\l}$.
\end{proof}
\subsection{The periodic Full Kostant-Toda lattice is a  Hamiltonian system}
 We intend in this section to show that the periodic Full 
Kostant-Toda is a Hamiltonian system for this Poisson structure.
 But, a small difficulty appears here: the function on $\FF(\tilde{\g}) $ with respect 
to which this equation is the Hamiltonian:
\begin{equation}\label{hamifktp}
H(L(\lambda)):=\frac12{\inn{L(\lambda)}{L(\lambda)}}_{\l},
\end{equation}%
which is not an element of $\FF(\tilde{\g})$. Fortunately, there exist an elements of $F_H \in \FF(\tilde{\g})$ 
whose restriction
to $\TT_\l$ is equal to the restriction of $H$, for instance the function
\begin{equation}\label{hamifktpFake}
 F_H (x(\l)):=\frac12 \left( \inn{x_{-1}}{x_1}+\inn{x_0}{x_0}+\inn{x_1}{x_{-1}} \right),
\end{equation}%
where $x(\l)= \sum_{i \in \Z}x_i \l^i$.
 We define the \emph{Hamiltonian vector fields of $H$ on $\TT_\l$} 
(or of any function on $\tilde{g}$ which satisfies the same property) to be the  Hamiltonian vector 
field (on $\TT_\l$) of any of these functions
 (Hamiltonian vector field which does not depend of the choice of $F_H$, since by Proposition \ref{fgyu} 
the Hamiltonian vector field of a function that vanishes on $\tilde \g$ also vanishes on $\TT_\l$).
\begin{proposition}%
 The Hamiltonian vector field  on $\TT_{\l}$  of the function $H$ defined in (\ref{hamifktp}) coincides with the equation
 of motion (\ref{mvpf}) of the periodic  Full Kostant-Toda lattice.
\end{proposition}
\begin{proof}
This proposition is just a particular case of the Adler-Kostant-Symes theorem \cite[Theorem 4.37]{Liv},
 up to the fact that we have to adapt it to the infinite dimensional setting.
 By definition, the Hamiltonian vector field on $\TT_{\l}$ of the function $H$ 
is the Hamiltonian vector field of the function $ F^H $ introduced in (\ref{hamifktpFake}).
Since the gradient of $ F^H ( x(\l) ) $ at a point $x(\l) \in \tilde{\g}$ is $ x_{-1}\l^{-1} + x_0 + x_1 \l$,
we have $ \nabla_{L(\l)} F^H = L(\l)  $ for every $L(\l) \in \TT_\l \subset \g \l^{-1} + \g + \g \l$, so that
 $$ \X_{ H}(L(\l)) = \frac{1}{2}\lb{ \tilde{R}(L(\l)),L(\l)} = \frac{1}{2} \lb{L(\l)_+ - L(\l)_-,L(\l)},$$
by definition of $\tilde{R}$. Hence
 \begin{eqn*}
\X_H(L(\lambda))=-\lb{(L(\lambda))_-,L(\lambda)}.
\end{eqn*}%
\end{proof}
\section{The Liouville integrability of the periodic Full Kostant-Toda lattice}
\label{sect3}
As in  Section \ref{sect1}, we choose  $\g$  a simple Lie  algebra, equipped with the Killing form  $\INN$, and
  $\h$  a Cartan subalgebra. Let  $P_1,\dots, P_{\ell}$ be  a generating family of the algebra of the  $\ad$-invariant
 polynomial
functions on $\g$, such that the degree of $P_i$ is $m_i+1$, for all $1\leqslant i\leqslant \ell$,
 where $m_1,\dots,m_{\ell}$ are
the exponents of $\g$ (we notice that $m_1\leqslant \dots\leqslant m_{\ell}$).
Each ${ P}_i$ extends  on $\tilde \g$ to a function
 ${\tilde P}_i$  with values in  $\C[\l,\l^{-1}]$, each of these functions is an
 $\ad$-invariant function  of $\tilde \g$ with  values
 in  $\C[\l, \l^{-1}]$, so each  coefficient at  $\l$ is an $\ad$-invariant function on  $\tilde\g$
with  value in  $\C$. Let  $\tilde{F}_{j,i}$ be a functions on $\tilde{\g}$, defined by:
\begin{equation}\label{equinfi}
{\tilde P}_i(L(\lambda))=\sum_{j=-\infty}^{\infty}\lambda^{-j}{\tilde
  F}_{j,i}(L(\lambda)), \qquad  \forall L(\l)\in\tilde\g.
\end{equation}%
\begin{remark}\label{remark-involv0}
 Let  $H$ be the Hamiltonian of the periodic  Full Kostant-Toda lattice, defined in  (\ref{hamifktp}) by:
\begin{eqn*}
H(x(\lambda))=\frac12{\inn{x(\lambda)}{x(\lambda)}}_{\l}, \qquad\qquad \forall x(\l)\in\tilde \g.
\end{eqn*}%
It is clear that  $H$ is  homogeneous, $\ad$-invariant of degree $2=m_1+1$, therefore  we can take   ${\tilde P}_1:=H$.
\end{remark}
 The functions ${\tilde F}_{j,i}$, for $1\leqslant i\leqslant \ell$ and $j\in\Z$, are $\ad$-invariant functions on $\tilde\g$.
 According to the AKS Theorem \cite[Theorem 4.36]{Liv}, they should in involution for the 
 Poisson $\tilde R$-bracket~${\PB}_{\tilde R}$. However, there is a technical issue here: 
strictly speaking, one can not apply the AKS theorem, since our Lie algebra is infinite dimensional
 and, moreover, the functions ${\tilde F}_{j,i}$ are not in $\FF(\tilde{\g})$ in general.
 The conclusion the AKS theorem, however, holds, at least after restriction to $\TT_\l$.
\begin{proposition} \label{prop:involv}
The restrictions to $\TT_\l$ of the functions $({\tilde F}_{j,i})$, $1\leqslant i\leqslant \ell$, $j\in\Z$, pairwise commute.
\end{proposition}
 \begin{proof}
The proof is an adaptation of the proof of the AKS theorem. For all $1\leqslant i\leqslant \ell$, $j\in\Z$,
 there exists a function $F^{{\tilde F}_{j,i}} \in \FF(\tilde{g})$ such that $F^{{\tilde F}_{j,i}}  $ and  ${\tilde F}_{j,i}$
 coincide on $ \TT_\l$.
 Moreover, although $F^{{\tilde F}_{j,i}}$ is not $\ad$-invariant on $\tilde \g$, we can assume that at
 all point $x(\l) \in \TT_\l$:
 \begin{equation}\label{eq:fake_ad_inv}
 \lb{ x(\l) , \nabla_{x(\l)} F^{{\tilde F}_{j,i}} } = 0.
\end{equation}%
For instance, the function ${\tilde F}_{j,i} \circ p_n$, where  $p_n$ is the projection of
 $\tilde{\g}$ on $ \sum_{i=-n}^n \l^i \g$,
satisfies these conditions for $n$ large enough.\\
 Since for all possible indices $ F^{{\tilde F}_{j,i}}$ and ${\tilde F}_{j,i}$ coincide when
 restricted to the Poisson submanifold $\TT_\l$, the Poisson brackets $\pb{ {\tilde F}_{j,i}, {\tilde F}_{k,l}}_{\tilde{R}}$
 and $\pb{ F^{{\tilde F}_{j,i}}, F^{{\tilde F}_{k,l}}}_{\tilde{R}}$ coincide on $\TT_\l$
for all possible indices, so that we are left with the task of proving that 
$\pb{ F^{{\tilde F}_{j,i}}, F^{{\tilde F}_{k,l}}}_{\tilde{R}}=0$ on $\TT_\l$. From now, 
the usual computation that proves of AKS theorem \cite[Theorem 4.36]{Liv} can be repeated word by word:
 \begin{eqnarray*} 
 \renewcommand{\arraystretch}{0.9}
\pb{ F^{{\tilde F}_{j,i}}, F^{{\tilde F}_{k,l}}}_{\tilde{R}} (x(\l))
&=&  \inn{x(\l)}{ \lb{ \nabla_{x(\l)} F^{{\tilde F}_{j,i}} ,\nabla_{x(\l)} F^{{\tilde F}_{k,l}} }_{\tilde{R}} }_\l  \\
 &=&  \frac12 \inn{x(\l)}{ \lb{ {\tilde{R}} (\nabla_{x(\l)} F^{{\tilde F}_{j,i}}) ,\nabla_{x(\l)} F^{{\tilde F}_{k,l}} } }_\l
 \\ & &+ \frac12 \inn{x(\l)}{ \lb{ \nabla_{x(\l)} F^{{\tilde F}_{j,i}} , {\tilde{R}}(\nabla_{x(\l)} F^{{\tilde F}_{k,l}} )} }_\l \\
 &=&  -\frac12 \inn{  \lb{ x(\l),\nabla_{x(\l)} F^{{\tilde F}_{k,l}} } }{  {\tilde{R}} (\nabla_{x(\l)} F^{{\tilde F}_{j,i}})  }_\l  \\
 & &+ \frac12 \inn{ \lb { x (\l),  \nabla_{x(\l)} F^{{\tilde F}_{j,i}}} }{ {\tilde{R}}(\nabla_{x(\l)} F^{{\tilde F}_{k,l}} )}_\l \\
 &=& 0
    \end{eqnarray*}%
where, in the last line, we have used twice (\ref{eq:fake_ad_inv}).
  \end{proof}
\begin{remark}\label{remark-involv2}
 There is therefore a large number of functions in involution  that are a goods candidates for
 the integrability of the periodic Full Kostant-Toda lattice. It will be show later that most of them are zero or constants and the remaining functions give the exact integrability.
\end{remark}
In this section we some results that we give in the following lemma.
\begin{lemma}\label{idensum}
Let  $\g$ be a simple Lie algebra  of rank $\ell$,  $\h$ be  a Cartan subalgebra of
 $\g$, $\Phi$ be a system of roots of   $\g$ associated to  $\h$, $(\a_1,\dots,\a_{\ell})$ be  a  basis of   $\Phi$
and   $h_1,\dots,h_\ell$ be the corresponding to simple  coroots. For every  $\gamma\in\Phi$, we choose
  $e_{\gamma}$  a non-zero eigenvector of~$\gamma$. Let
\begin{eqn*}
(x_1,\dots,x_{\ell})\cup (x_{\gamma})_{  \gamma\in\Phi}
\end{eqn*}%
be the coordinates system on  $\g$ given, for every  $ 1\leqslant i\leqslant
\ell$ and every  $ \gamma\in\Phi$  and for every  $x\in\g$, by:
\begin{eqn*}
\left\{
\renewcommand{\arraystretch}{0.9}
\begin{array}{l}
x_i(x)=\inn{h_i}{x},\\
x_{\gamma}(x)=\inn{e_{-\gamma}}{x}.
\end{array}
\right.
\end{eqn*}%
Let  $P$ a homogeneous $\ad$-invariant polynomial on  $\g$ of  degree $m+1$.\\
\noindent
\emph{(1)} The polynomial  $P$ is a linear combination of the monomials of the following form
\begin{equation}\label{kopr23}
x_{\gamma_1}\dots x_{\gamma_k}x_{p_1}\dots x_{p_{j}},
\end{equation}%
where $p_1,\dots,p_j\in\{1,\dots,\ell\}$  and such that:
\begin{equation}\label{twdx1}
\left\{
\begin{array}{c}
k+j=m+1,\\
\sum_{i=1}^{k}|\gamma_i|=0.
\end{array}
\right.
\end{equation}%
\emph{(2)} Let
 $i_1,\dots,i_p \in \{1,\dots,\ell\}$ and
  $\gamma_1,\dots,\gamma_{q}\in\Phi$, where  $p,q\in\N$.
 If
\begin{equation}\label{eq:conditions}
m+1-(p+q)+\sum_{i=1}^{q}|\gamma_i|< 0\qquad \textrm{ or } \qquad
\sum_{i=1}^{q}|\gamma_i|> 0
\end{equation}%
then, for every   $y\in\h\oplus \g_1$,
\begin{equation}\label{eq:voulue}
\can{\diff_y^{p+q}P,(h_{i_1},\dots,h_{i_{p}},e_{\gamma_1},\dots,e_{\gamma_{q}})}=0.
\end{equation}
\end{lemma}
\begin{proof}
{\em(1)} Every  homogeneous polynomial  of degree
 $m+1$ is a linear combination of monomials of the form
 (\ref{kopr23}) with   $k+j = m+1$. We need  to show that  when this polynomial is  $\ad$-invariant, the second
 condition  of  system (\ref{twdx1}) is satisfied for every  monomial that  appear in its decomposition.\\
Let
$h\in\h$ be such that  $\a_i(h)=1$ for every  $i=1,\dots,\ell$.  We define  a linear vector field
 $\widetilde{\ad}_h$ on  $\g$ by:
$$
\widetilde{\ad}_h [F] (x):=\can{\diff_x F,\ad_hx}=\inn{\nabla_x F}{\ad_h x},
$$
for every  $F \in \FF (\g)$ and every   $x\in \g$. On the one hand,  for every  $\gamma\in\Phi$
\begin{eqn*}
\widetilde{\ad}_h[x_{\gamma}](x)=\inn{\ad_hx}{e_{-\gamma}}
                       =\gamma(h)x_{\gamma}(x)
                       =|\gamma|\, x_{\gamma}(x)
\end{eqn*}%
while   $\widetilde{\ad}_h[x_i]=\inn{\ad_hx}{h_i}=0$, for
$i\in\{1,\dots,\ell\}$ on the other hand. These two  properties  imply 
\begin{eqnarray}
\renewcommand{\arraystretch}{0.9}
\widetilde{ \ad}_h[x_{\gamma_1}\dots x_{\gamma_k}x_{p_1}\dots x_{p_{j}}]&=&\sum_{i=1}^{k}\tilde{ \ad}_h[x_{\gamma_i}]x_{\gamma_1}\dots \hat{x}_{\gamma_i}\dots x_{\gamma_k}x_{p_1}\dots x_{p_{j}}\nonumber\\
                 &=&(\sum_{i=1}^{k}|\gamma_i|)x_{\gamma_1}\dots x_{\gamma_k}x_{p_1}\dots x_{p_{j}}.\label{edentop}
\end{eqnarray}%
 Since  $P$ is  an $\ad$-invariant polynomial,
 $\widetilde{\ad}_h [P](x)=\inn{\ad_h x}{\nabla_x P}=\inn{h }{[x,\nabla_x P]}=0$.
 Therefore, according  to   (\ref{edentop}), the  sum $\sum_{i=1}^{k}|\gamma_i|$
 vanishes for each monomial appearing in the decomposition of  $P$.\\
{\em(2)}
 If  $p+q \geq m+2$,  equation  (\ref{eq:voulue}) holds   automatically, because the  degree
 of  $P$ is  $m+1$.
We assume  for   $p+q \leqslant m+1$, the first  point of the  lemma  implies that, for every  $y\in \g$  and  every
 homogeneous elements $z_1,\dots,z_{m+1} \in \g$
 with
$  \sum_{k=1}^{m+1} |z_i| \neq 0$,
 \begin{equation}\label{eq:equiv_homog}
 \can{\diff_y^{m+1}P,(z_{1},\dots,z_{m+1} )}=0.
\end{equation}%
Let  $i_1,\dots,i_p\in\{1,\dots,\ell\}$ and let  $\gamma_1,\dots,\gamma_q\in\Phi$. Since the function  $$y \mapsto
\can{\diff_y^{p+q}P,(h_{i_1},\dots,h_{i_{p}},e_{\gamma_1},\dots,e_{\gamma_{q}})},$$
is homogeneous of  degree  $m+1-p-q$, according to Taylor   formula, it is equal to
$$ y \mapsto \frac{1}{(m+1-p-q)!}\can{\diff_y^{m+1}P,(h_{i_1},\dots,h_{i_{p}},e_{\gamma_1},\dots,e_{\gamma_{q}},y^{m+1-p-q})}.$$
By restricting to  $\h \oplus \g_1$, this last function is a linear combination  of monomials of the form
  $$ x_1^{a_1} \dots,x_\ell^{a_\ell} x_{\alpha_1}^{b_1} \dots
  x_{\alpha_\ell}^{b_\ell}  $$
where  $\sum_{k=1}^\ell (a_k +b_k)= m+1-p-q $.  The coefficient in the decomposition of $P$  of the above
 monomial is
 $$   \frac{1}{(m+1-p-q)!} \can{\diff_y^{m+1}
P,(h_{i_1},\dots,h_{i_{p}}, h_1^{a_1},\dots, h_\ell^{a_\ell},e_{\gamma_1},\dots,e_{\gamma_{q}},
  e_{\alpha_1}^{b_1},\dots,e_{\alpha_\ell}^{b_\ell} )}.$$
According to  (\ref{eq:equiv_homog}), this coefficient vanishes if
 $$ \sum_{i=1}^q |\gamma_i|+ \sum_{k=1}^{\ell} b_k \neq 0. $$
Since  $ \sum_{k=1}^{\ell} b_k \in \{0,\dots,m+1-p-q \} $, all the  coefficients vanish if one of the two
 conditions   (\ref{eq:conditions}) is satisfied.
  \end{proof}
\begin{proposition}\label{biyt4}
For  $i=1,\dots,\ell,$
the  restriction   of   ${\tilde P}_i$ to
$\TT_{\l}$ is given by
\begin{equation}\label{defIS}
{\tilde P}_i(L(\lambda))=\sum_{j=0}^{m_i}\lambda^{-j}{\tilde
  F}_{j,i}(L(\lambda))+\l \, c\,\delta_{i,\ell},\qquad\qquad  \forall L(\l)\in\TT(\l),
\end{equation}%
where
$c$ is  a non-zero  constant.
\end{proposition}
\begin{proof}
Since the degree of ${\tilde P}_i$, for all $1\leqslant i\leqslant \ell$ is
equal  to  $m_i+1$,
 the restrictions of the functions $\tilde{F}_{k,i}(L(\l))$  (constructed in  (\ref{equinfi})) to  $\TT_{\l}$  vanish for every $1\leqslant i\leqslant \ell$ and every
 $-m_i-1\leqslant j\leqslant m_i+1$ and
\begin{eqn*}
{\tilde P}_i(L(\lambda))=\sum_{k=-m_i-1}^{m_i+1}\lambda^{-k}{\tilde
  F}_{k,i}(L(\lambda)).
\end{eqn*}%
Let show that   $\tilde {F}_{m_i+1,i}$  vanish on  $\TT_{\l}$,  for every  $1\leqslant i \leqslant\ell$. Let
 $L(\l)=\l e_{-\b}+X+\l^{-1}Y\in\TT_{\l}$, we notice that
\begin{eqn*}
\tilde {P_i}(L(\l ))=\l^{-m_i-1}\tilde{P}_{i}(Y+\l^2e_{-\b}+\l X ).
\end{eqn*}%
Therefore the coefficient of degree $-m_i-1$ is
\begin{eqn*}
{\tilde F}_{m_i+1,i}(L(\l ))=P_i(Y).
\end{eqn*}%
Since   $Y$ is  an  element of  $\g_{>0}$, it is   nilpotent. This implies,  according to   \cite[Theorem 8.1.3]{Dixmier}
 that  $P(Y)$ is zero
for every    $P$ an  $\Ad$-invariant polynomial on  $\g$.
 \smallskip

Let us show that the  functions  $\tilde{F}_{j,i}$,  for all   $j$
 strictly lower to   $-1$  vanish  and that the  function  $\tilde
{F}_{-1,i}$ vanish except for  $i=\ell$, in which case it is a constant function.\\ 
The  extensions  ${\tilde x}_i$ and  ${\tilde x}_{\gamma}$, for every  $1\leqslant i\leqslant \ell$ and every
 $\gamma\in\Phi$ to  $\tilde \g $, of the  coordinate functions $(x_i, x_{\gamma}, 1\leqslant i\leqslant \ell, \gamma\in\Phi)$ on
  $\g$ defined in    Lemma  \ref{idensum}  have  restrictions to  $\TT_{\l}$ given by:
\begin{equation}\label{systcoo}
\renewcommand{\arraystretch}{1.5}
\left\{
\renewcommand{\arraystretch}{0.9}
\begin{array}{lll}
x_i,& 1\leqslant i\leqslant \ell, & (\textrm{type } I)\\
x_{{-\gamma}},&\textrm{if }\gamma\in\Phi_+\backslash{\b}, & (\textrm{type } II) \\
 x_{-\beta}+\l, &\textrm{if }  \gamma=\beta, & (\textrm{type } III)\\
\l^{-1}y_{\gamma}+1, &\textrm{if } \gamma\in\Pi, & (\textrm{type } IV)\\
\l^{-1}y_{\gamma},&\textrm{if }\gamma \in\Phi_+\backslash \Pi, & (\textrm{type } V)
\end{array}
\right.
\end{equation}%
here $y_{\gamma}$ stands for  $x_{\gamma}$ for any $\gamma$ a positive  root.
Then, for each  $P_i$  an  $\Ad$-invariant homogeneous polynomial on  $\g$ of degree $m_i+1$, the
 restriction to  $\TT_{\l}$  of  its  extension ${\tilde P}_i$ on  $\tilde \g$ is a combination of monomials of the
following form
\begin{equation}\label{sffrtt}
\renewcommand{\arraystretch}{0.9}
\begin{array}{l}
x_{p_1}\dots x_{p_{h}},\\
\times\\
x_{-{\gamma_1}}\dots x_{-\gamma_p},\\
\times\\
( x_{-\b}+\l)^l, \\
\times\\
(\l^{-1}y_{\a_{j_1}}+1)\dots(\l^{-1}y_{\a_{j_k}}+1)\\
\times\\
\l^{-1}y_{\delta_1}\dots\l^{-1}y_{\delta_q},
\end{array}
\end{equation}%
where  $\a_{j_1},\dots\a_{j_k}\in\Pi $,
$\gamma_1,\dots,\gamma_p\in\Phi_+\backslash\b$,
   $ \delta_1,\dots,\delta_{q}\in\Phi_+\backslash \Pi$,
 $l\in\N$ et $p_1\dots,p_h\in\{1,\dots,\ell\}$ and  where the following conditions are satisfied:
\begin{eqn*}
\left\{
\begin{array}{ll}
h+p+l+k+q= m_i+1 & \textrm{ (C1)}, \\
-\sum_{i=1}^{p}|\gamma_i|-l|\b|+k+\sum_{i=1}^{q}|\delta_i|=0 &\textrm{ (C2) }.
\end{array}
\right.
\end{eqn*}%
Of course, it should be understood  that if  $h=0$ or
  $p=0$ or $j=0$ or $k=0$ or $q=0$, then in (\ref{sffrtt}) the corresponding term is equal to $1$.\\
The first condition simply comes from the fact that $P_i$ is homogeneous of degree $m_i+1$ and the second is
 a consequence of the first point of Lemma \ref{idensum}, claiming that the $P_i$ are homogeneous
 of degree zero with respect to  the root weight.\\
Let us now show that the  functions  $\tilde{F}_{j,i}$ vanish, for every   $j$ strictly lower to  $-1$. For
 all  $1\leqslant  i\leqslant q $ the length of the root $\delta_i$ is lower than or equal to $|\b|=m_{\ell}$.
 Furthermore, $k$ is lower than or equal to  $m_{i}+1$, hence to $m_{\ell}+1$. But we can not have  $k=m_{\ell}+1$,
because that   implies $h=p=l=q=0$ and  contradicts the  second condition~(C2).
 Therefore  $k \leqslant m_{\ell}$,  and we obtain
 the inequality
\begin{equation}
\label{eq:qetj}
S=\sum_{i=1}^{p}|\gamma_i|= -lm_{\ell}+k+\sum_{i=1}^{q}|\delta_i|\leqslant (1+q-l)m_{\ell}.
\end{equation}%
The length of roots   $\gamma_1,\dots,\gamma_p$  is
positives, their  sum  $S$ is   positive (or  zero when   $p=0$). Hence   $l\leqslant q+1$.
This  implies that  the monomials that make up the restriction to $\TT_\l $ of $\tilde{P}_i $ have at least $l-1$
products of functions of type V whenever they have  $l$ products of the functions of type III. This product contains
  one and only one a term in $\lambda^j$ for $j\geq 1 $. Since the   other types (I-II-IV)
are polynomials in $ \l ^{-1} $, the restriction to $ \TT_\l $ of $ \tilde{P}_i $
 contains only  a term in $\lambda^j $ for $ j\geq 1 $, i.e., the restriction of  the functions  $\tilde {F}_{j,i}$ vanish
for every  $j \leqslant -2$.\smallskip\\
We now show that the function $\tilde {F}_ {-1,i} $ vanish except for $ i = \ell $ in which case it is a non-zero constant.
It follows from  (\ref{systcoo})
that a term in  $ \lambda $ appears in the monomials which
compose
 $ {\tilde {P}}_i $ that  if $ l \geq q +1 $. But we now that  $l \leqslant q+1$, then  $l=q+1$. According to  (\ref{eq:qetj}),
this implies that  $p=0$, and  that  $j=m_{\ell}$. Hence the condition (C1) becomes $ h+2q+1+m_{\ell}=m_i+1 $, this
 in turn implies $m_i=m_{\ell}$ and  $h=q=0 $, then  $l=1$. The monomials where the term  in $\lambda$ 
appeares are therefore  the  product
  of $ m_{\ell} $ terms of   the type IV with one term  of the type III, i.e.,  is the  product
\begin{eqn*}
(x_{-\b}+\l )(\l^{-1}y_{\a_{j_1}}+1)\dots(\l^{-1}y_{\a_{j_{m_{\ell}}}}+1),
\end{eqn*}%
where  $\alpha_{j_1},\dots,\alpha_{j_{m_{\ell}}}$ are a simple roots. But the coefficient  in  $ \l $ appearing
in this case is
constant.
\end{proof}
Most of the functions ${\tilde F}_{j,i}, 1\leqslant  i\leqslant \ell, j\in\Z $ are identically zero (or constant) after restriction  to $\TT_{\l}$. For the remaining functions, we introduce the following notation.\\
{\bf{Notation}:}
 We denote  by $\tilde{\FF}_{\l}$ the   family
of the restriction of  functions  $\tilde{ F}_{j,i}$ to  $\TT_{\l}$, for every  $1\leqslant
  i\leqslant\ell$ and every  $0\leqslant j\leqslant m_i$, i.e.,
 \begin{equation}\label{fam-integrable}
{\tilde \FF}_{\lambda}:=({\tilde F}_{j,i},~~ 1\leqslant i\leqslant \ell,~~ 0\leqslant j\leqslant m_i).
\end{equation}%
We can now give the main result of this article.
\begin{theorem}\label{principal-theorem}
The triplet $(\TT_{\l},{\tilde\FF}_{\lambda},{\PB}_{\tilde R})$ is an integrable system.
\end{theorem}
\begin{proof}
According to the definition of integrability in the sense of Liouville (see \cite[Definition 4.13]{Liv})
to prove Theorem (\ref{principal-theorem}), we must show that:
\begin{enumerate}
\item[(1)] ${\tilde \FF}_{\l}$ is involutive for the Poisson ${\tilde R}$-bracket ${\PB}_{\tilde R}$.
\item[(2)] ${\tilde \FF}_{\l}$ is independent on $\TT_{\l}$.
\item[(3)] The cardinal of ${\tilde \FF}_{\l}$ satisfies
\begin{equation}
\card {\tilde\FF}_{\l}=\dim\TT_{\l}-\frac 12 \Rk(\TT_{\l},{\PB}_{\tilde\R }).
\end{equation}%
\end{enumerate}
The proofs of these three points  are given in respectively  Proposition \ref{prop:involv},
Proposition \ref{funinde} and Proposition \ref{prkhim1}, the latter two propositions being given in the next two subsections.
\end{proof}%
\subsection{The family ${\tilde \FF}_{\l}$ is independent on $\TT_{\l}$}
We use  an unpublished result  of Ra\"is \cite{Rais}, which  establishes  the independence
 of a large family of functions on $\g\times\g$.  We  stated  this result below and refer to  \cite[Section 1]{ddg} for
a proof.
\begin{theorem}\label{indvki}
Let  $P_1,\dots,P_{\ell}$   be a generating family of homogeneous polynomials of the algebra of  $\Ad$-invariant polynomial functions on   $\g$.
 Let  $e$ and  $h$ be  two elements of $\g$,
 such that  $e$ is regular and  $[h,e]=2e$.\\
 For every  $F\in\FF(\g)$, and every  $y\in\g$, we denote  by  $\diff^k_yF$
 the differential of  order  $k$ of $F$  at $y$.
 Denote by  $V_{k,i}$, for every   $1\leqslant i\leqslant \ell$ and  $0\leqslant k\leqslant
 m_i$,  the elements of  $\g$ defined by:
\begin{equation}\label{dpvki}
\inn{V_{k,i}}{z}=\can { \diff^{k+1}_hP_i, (e^k,z)},\qquad \qquad \forall z\in\g,
\end{equation}%
where, for every $x\in\g$ and $k\in\N$, $x^k$ is a shorthand for $(x,\dots,x)$ ($k$ times).
\\
\emph{(1)} The family   $\FF_1:=(V_{k,i}, 1\leqslant i\leqslant \ell$ and $ 0\leqslant k\leqslant m_i)$
is  linearly   independent;\\
\emph{(2)}
The subspace  generated by    $\FF_1$ is the Lie  subalgebra formed by the sum of the all  eigenspaces of $ \ad_h$  associated with positive or zero eigenvalues.
\end{theorem}
We now   show the independence of the differentials of the family of functions $\tilde \FF_{\l}$ defined in  \ref{fam-integrable}  in a particular
  point  of  $\TT_{\l}$ (which implies  the independence of the family ${\tilde
   \FF}_{\ell}$  because its elements are  polynomials).
\begin{proposition}\label{funinde}
The family of functions $\tilde\FF_{\lambda}$ is  independent on $\TT_{\l}$.
\end{proposition}
\begin{proof}
Let   $h\in\h$, such that   $[h,e]=2e$. We first  prove that  $\tilde
{\FF}_{\l}$ is  independent  at the point  $L_1(\l):=\l e_{-\beta}+ h+e+\l^{-1} e$.\\
 We compute the differential of the function  ${\tilde P}_i$ (valued in  $\C[\l,\l^{-1}]$) at the  point $L_1(\l)$.
 Let
 $ a(\l):=A+\l^{-1}B\in T_{L_1(\l)}\TT_{\l}=\bigoplus_{-|\b|\leqslant i\leqslant
   0}^{}{\tilde\g}_i$, we have the equality:
\begin{eqnarray}
\renewcommand{\arraystretch}{0.9}
\can{\diff_{L_1(\l)}{\tilde
    P}_i,a(\l)}&=&\can{\diff_{h+(1+\l^{-1})e+\l e_{-\b}}{\tilde
    P}_i,a(\l)}\nonumber\\
            &=&\sum_{j=0}^{m_i}\frac{\l^j}{j!}\can{\diff_{h+(1+\l^{-1})e}^{j+1}{\tilde
      P}_i,((e_{-\b})^j,a(\l ) )}\nonumber\\
           &=&\can{\diff_{h+(1+\l^{-1})e}{\tilde
      P}_i,a(\l)}+\sum_{j=1}^{m_i}\frac{\l^j}{j!}\can{\diff_{h+(1+\l^{-1})e}^{j+1}{\tilde  P}_i,((e_{-\b})^j,a(\l ))}\nonumber\\
         &=&\can{\diff_{h+(1+\l^{-1})e}{\tilde  P}_i,A}+\l^{-1}\can{\diff_{h+(1+\l^{-1})e}{\tilde P}_i,B}\nonumber\\
       &&+\sum_{j=1}^{m_i}\frac{\l^j}{j!}\can{\diff_{h+(1+\l^{-1})e}^{j+1}{\tilde
      P}_i,((e_{-\b})^j,A)}\nonumber \\&&+\sum_{j=1}^{m_i}\frac{\l^{j-1}}{j!}\can{\diff_{h+(1+\l^{-1})e}^{j+1}{\tilde P}_i,
((e_{-\b})^j,B)}.\label{eqsdet}
\end{eqnarray}%
To go from the first to the  second line, we  have used the  fact that the polynomial $\tilde{P}_i$ has degree $m_i+1$
(therefore its differential is of degree $m_i$).\\
Since $A\in\g_{\leqslant 0}$, it is of the form  $A=\sum_{i=1}^{\ell}a_ih_i+\sum_{\gamma\in\Phi_+}^{}a_{\gamma}e_{-\gamma}$. Since for
$1\leqslant j\leqslant m_i$ the integers, respectively $m_i+1-j-1+j|-\b|+|h_i|$ and  $m_i+1-j-1+j|-\b|+|e_{-\gamma}|$, which are
smaller or equal, respectively to  $ -j -m_{\ell}(j-1) $ and  $ -j -m_{\ell}(j-1) +|e_{-\gamma}|$ are strictly negatives.
  According to the second item of Lemma   \ref{idensum}, therefore: 
\begin{equation}\label{eqfpt1}
\sum_{j=1}^{m_i}\frac{\l^j}{j!}\can{\diff_{h+(1+\l^{-1})e}^{j+1}{\tilde
      P}_i,((e_{-\b})^j,A)}=0.
\end{equation}%
Moreover,  $B\in\g_{>0}$ is   of the  form    $B=\sum_{\gamma\in\Phi_+}^{}b_{\gamma}e_{\gamma}$.
 By using  again the second item of Lemma~\ref{idensum}, we deduce that:
\begin{equation}\label{eqfpt2}
\can{\diff_{h+(1+\l^{-1})e}{\tilde P}_i,B}=0.
\end{equation}%
Using  Equations  (\ref{eqfpt1}) and  (\ref{eqfpt2}),   (\ref{eqsdet}) becomes:
\begin{equation}\label{eqsdet1}
 \can{\diff_{L_{1}(\l)}{\tilde
    P}_i,a(\l)}=\can{\diff_{h+(1+\l^{-1})e}{\tilde P}_i,A}+
\can{\sum_{j=1}^{m_i}\frac{\l^{j-1}}{j!}\diff_{h+(1+\l^{-1})e}^{j+1}{\tilde
      P}_i,((e_{-\b})^j,B)}.
\end{equation}%
We denote by  ${\tilde H}_{j,i}$ the  function  defined on  $\g\times\g$  by:
\begin{eqn*}
{\tilde P}_i(X+\l^{-1}Y)=\sum_{j=0}^{m_i+1}\l^{-j}{\tilde H}_{j,i}(X,Y), \qquad\qquad \forall X,Y\in\g\times\g.
\end{eqn*}%
We clearly have:
\begin{eqn*}
{\tilde P}_i(X+(1+\l^{-1})Y )=\sum_{j=0}^{m_i+1}\l^{-j}{\tilde H}_{j,i}(X+Y,Y).
\end{eqn*}%
We  notice that on  $\g\times\g_{>0}$,
\begin{enumerate}
\item[(1)]  The  function ${\tilde H}_{m_i+1,i}(X+Y,Y)=P_i(Y)=0$;
\item[(2)]  The  differentials of  ${\tilde H}_{0,i},\dots,{\tilde H}_{m_i,i}$
 at  point $(h+e,e)$ do not depend on the variable $Y$, because according to (\ref{eqfpt2}),
  $\can{\diff_{h+(1+\l^{-1})e}{\tilde P}_i,B}=0, \forall B\in\g_{>0}$.
\end{enumerate}
Theses two   points implies that
\begin{equation}\label{pfre1}
\diff_{h+(1+\l^{-1})e}{\tilde P}_i=\sum_{j=0}^{m_i}\l^{-j}\pp{{\tilde H}_{j,i}}{X}{(h+e,e)},
\end{equation}%
where  $\pp{{\tilde H}_{j,i}}{X}$, for every  $1\leqslant i\leqslant\ell $ and  $0\leqslant
j\leqslant m_i$, stands for  the  differential  of  $\tilde{H}_{j,i}$ with respect to the first  variable.
 Using  Equation  (\ref{pfre1}),  Equation (\ref{eqsdet1}) becomes:
 \begin{eqnarray}\label{pfre4}
 \renewcommand{\arraystretch}{0.9}
 \can{\diff _{L_1(\lambda)}{\tilde P}_i,
 a(\l)}&=&\sum_{j=0}^{m_i}\l^{-j}\can{\pp{{\tilde
     H}_{j,i}}{X}{(h+e,e)}, A}\nonumber\\
&&+\sum_{j=1}^{m_i}\frac{\l^{j-1}}{j!}\can{\diff_{h+(1+\l^{-1})e}^{j+1}{\tilde
      P}_i,((e_{-\b})^j,B)}.
 \end{eqnarray}%
 Since  $L_1(\l)$ is an  element of $\TT_{\l}$, according to  Relation   (\ref{defIS}),
 \begin{equation}\label{pfre3}
  \diff _{L_1(\lambda)}{\tilde
    P}_i=\sum_{j=0}^{m_i}\l^{-j}\diff_{L_1(\l)}{\tilde F}_{j,i}.
 \end{equation}%
 By using  Equations (\ref{pfre4}) and  (\ref{pfre3}), we  conclude that
 \begin{eqnarray}\label{pfre2}
 \renewcommand{\arraystretch}{0.9}
 &&\sum_{j=0}^{m_i}\l^{-j}\can{\diff_{L_1(\l)}{\tilde
     F}_{j,i},a(\l)}=
\sum_{j=0}^{m_i}\l^{-j}\can{\pp{{\tilde H}_{j,i}}{X}(h+e,e), A}\nonumber \\
&&+\sum_{j=1}^{m_i}\frac{\l^{j-1}}{j!}\can{\diff_{h+(1+\l^{-1})e}^{j+1}{\tilde
      P}_i,((e_{-\b})^j,B)}.
 \end{eqnarray}%
 It suffices therefore to prove that $\pp{{\tilde H}_{j,i}}{X}(h+e,e)$ are independent as linear forms on $\g_{\leqslant 0}$.\\
Let  $h'=h+e$, since  $e=\sum_{i=1}^{\ell}e_i$ is a regular   element of
 $\g$ and  $[h',e]=e$, according  to the first point of   Theorem \ref{indvki}
 the family of linear form on $\g$
  \begin{equation}\label{eq:partial-derivative}
\pp{{\tilde H}_{0,i}}{X}(h',e),\dots,\pp{{\tilde H}_{m_i,i}}{X}(h',e) \qquad 1\leqslant i\leqslant \ell,
\end{equation}%
 is  independent. These linear forms are given by the
 gradients
 $V_{k,i},$ for  $ 1\leqslant i \leqslant \ell$ and  $0\leqslant k\leqslant m_i$,
 that belong to the
 space $E$ spanned by the eigenspaces
 of positive eigenvalues of $\ad_{h'} $
(see the second  point of
 Theorem \ref{indvki}). But the space spanned by the eigenspace of positive eignvalues of
both $\ad_h$ and $\ad_{h'}$ coincide with   $\g_{\geq 0} $.
Therefore the restrictions to $ \g_{\leqslant 0} $ of  the family (\ref{eq:partial-derivative})
  remain independent.
As a result, the differentials of the family of  functions  $ (\tilde F_{k, i}, 0 \leqslant i \leqslant m_i,
1 \leqslant i \leqslant \ell) $
 are independent at the  point $ L_1(\l) $ and therefore $ {\tilde \FF}_{\l} $ is independent on $ \TT_{\l} $.
\end{proof}
\subsection{The exact number on functions}
According to Equation (\ref{fam-integrable}), the cardinal of ${\tilde \FF}_{\l}$
 is related to the exponents $m_i$ of $\g$, $1\leqslant i\leqslant \ell$,
as  follows
\begin{equation}\label{rel-expo-card}
\card{\tilde \FF}_{\lambda}=\sum_{i=1}^{\ell}(m_i+1).
\end{equation}%
According to the classical relation  $\sum_{i=1}^{\ell}m_i=\frac{1}{2}(\dim\g-\ell)$ 
(see \cite[Theorem 7.3.8]{Dixmier}), Relation (\ref{rel-expo-card})
implies that   $\card{\tilde \FF}_{\lambda}=\frac{1}{2}(\dim\g+\ell)$. Moreover,  since the dimension
 of $\TT_{\l}$ is equal to $\dim\g$, the relation below is satisfied
\begin{eqn*}
\card {\tilde\FF}_{\l}=\dim\TT_{\l}-\frac 12 \Rk(\TT_{\l},{\PB}_{\tilde\R })
\end{eqn*}%
if and only if $\Rk(\TT_{\l},{\PB}_{\tilde\R })=\dim\g-\ell$. We need therefore to prove this last equality,
which shall be done in Proposition \ref{prkhim1} below.\\
~~\\
{\bf{The rank  of ${\PB}_{\tilde R}$ on $\TT_{\l}$}}\\
We show here  that there  exists $\ell$ independent Casimirs on $\TT_{\l}$ and  there exists a point
 $L_0(\l)$ of $\TT_{\l}$, such that the rank of the Poisson structure at this point is  $\dim\TT_{\l}-\ell=\dim\g-\ell$, which  proves that the rank of
the Poisson structure on $\TT_{\l}$ is  $\dim\g-\ell$.
\begin{proposition}\label{caspft}
The  functions ${\tilde F}_{m_1,1},\dots, {\tilde F}_{m_{\ell},\ell}$,  defined  in
(\ref{defIS}),  are
Casimirs  for the Poisson  $\tilde R$-bracket ${\PB}_{\tilde R}$.
\end{proposition}
We use Lemma \ref{lempft} below to show Proposition \ref{caspft}.
\begin{lemma}\label{lempft}
\emph{(1)} For every  $1\leqslant i\leqslant \ell$,  $Z(\l)=\sum_{k\geq
  0}\l^kZ_k\in\sum_{k\geq 0}^{}\l^k\g$  and    $
  Y\in\g_{>0}$, we have:
\begin{eqn*}
{\tilde F}_{m_i,i}(Z(\l)+\l^{-1}Y)=\can{\diff_Y{ P}_i,P_{\leqslant  0}(Z_0)},
\end{eqn*}%
where  $P_{\leqslant 0}$ is the projection of  $\g$ on  $\g_{\leqslant 0}$;\\
\emph{(2)}  At every   point   of  $\TT_{\lambda}$, the  gradients of the   functions
  $\tilde{F}_{m_1,1},\dots,\tilde{F}_{m_{\ell},\ell}$ are in ${\tilde \g}_+$.
\end{lemma}
\begin{proof}
{\em(1)} We denote,  for every  $k\in\N$ and every  $X(\l)\in\tilde\g$,  by $(X(\l))^k$ the  $k$-tuple
$(X(\l),\dots,X(\l))$ and for every  ${\tilde P}$, by  $\diff^k{\tilde P}_i$ the  $k^{\hbox{th}}$
differential  of ${\tilde P}_i$. The Taylor  formula of  $\tilde{P}_i$  at  point $Z(\l)+\l^{-1}Y$ is given by:
\begin{eqnarray}\label{dev12}
\renewcommand{\arraystretch}{0.9}
\tilde{P}_i(Z(\l )+\l^{-1}Y)&=&\l^{-m_i-1}\tilde{P}_i(\l Z(\l )+Y)\nonumber \\
                         &=&\sum_{j=0}^{m_i+1}\frac{\l^{j-m_i-1}}{j!}\can{\diff_{Y}^j{\tilde P}_i,(Z(\l ))^j}.
\end{eqnarray}
Recall from   (\ref{equinfi})  that the function $\tilde{F}_{m_i,i}$ is the  coefficient of degree  $-m_i$
 in  $\l$ of the  polynomial  $\tilde{P}_i$. Since   $Z(\l)\in\sum_{k\geq 0}\l^k\g$,  Formula  (\ref{dev12}) gives:
\begin{equation}
\tilde{F}_{m_i,i}(Z(\l)+\l^{-1}Y)=\can{\diff_Y\tilde{P}_i,Z_0}=\can{\diff_Y{P}_i,Z_0}.\label{eqyz}
\end{equation}%
  The  polynomial  $\can{\diff_Y{P}_i,Z_0}$ is  homogeneous of  degree $m_i+1$,
  of  degree $m_i$ with respect to the variable  $Y$ and of degree  $1$ with respect to the variable  $Z_0$.
  For all   $Y\in\g_{>0}$,    $\nabla_YP_i$ belong to   $\g_{\geq 0}$ hence:
\begin{eqn*}
\can{\diff_Y{P}_i,P_{> 0}(Z_0)}=0,%
\end{eqn*}%
where  $P_{>0}$ is the projection of  $\g$ on $\g_{>0}$. Therefore, Equation (\ref{eqyz}) becomes
\begin{eqn*}
\tilde{F}_{m_i,i}(Z(\l)+\l^{-1}Y)=\can{\diff_Y{P}_i,P_{\leqslant 0}(Z_0)},
\end{eqn*}%
where  $P_{\leqslant 0}$ is the projection of  $\g$ on  $\g_{\leqslant 0}$.\\
  {\em(2)} Let $X\in\g_{\leqslant 0}$,  $Y\in\g_{>0}$,   $L(\l)=\l
    e_{-\b}+X+e+\l^{-1}Y\in\TT_{\l}$ and let  $Z(\l)\in{\tilde \g}_{\geq 1}$. We recall that an
   element   $Z(\l)$ in  $\tilde{\g}_{\geq 1}$ has the following expression  $Z(\l) =\sum_{k\geq 0}\l^kZ_k
$, where  $Z_0\in\g_{\geq 1} $ and  $Z_k\in\g$ for all  $k>0$. According to the first  point  of the  lemma
\begin{eqn*}
\tilde{F}_{m_i,i}(L(\l))=\tilde{F}_{m_i,i}(L(\l)+Z(\l)), \qquad \forall Z(\l)\in\tilde{\g}_{\geq 1}.
\end{eqn*}%
The above equality  implies
\begin{eqn*}
{\inn{\nabla_{L(\l)}\tilde{F}_{m_i,i}}{Z(\l)}}_{\l}=0, \qquad \qquad \forall Z(\l)\in \tilde{\g}_{\geq 1}.
\end{eqn*}%
  This  implies that the  gradient of  $\tilde{F}_{m_i,i}$ at every  point of  $\TT_{\l}$ is in  $\tilde{\g}_+$.
\end{proof}
We now prove Proposition \ref{caspft}.
\begin{proof}
Let  $G\in\FF(\TT_{\lambda})$ and let
 $L(\lambda)\in\TT_{\lambda}$, we have:
\begin{eqnarray*}
\renewcommand{\arraystretch}{0.9}
{\pb{{\tilde F}_{m_i,i},G}}_{\tilde R}(L(\lambda ))&=&
{\inn{L(\l)}{{[\nabla_{L(\l )}{\tilde F}_{m_i,i},\nabla_{L(\l
          )}G]}_{\tilde R}}}_{\l}\\
          &=&{\inn{L(\l)}{[(\nabla_{L(\l )}{\tilde F}_{m_i,i})_+,(\nabla_{L(\l
          )}G)_+]-[(\nabla_{L(\l )}{\tilde F}_{m_i,i})_-,(\nabla_{L(\l
          )}G)_-]}}_{\l} \\
    &=&{\inn{L(\l)}{[\nabla_{L(\l )}{\tilde F}_{m_i,i},(\nabla_{L(\l
        )}G)_+]}}_{\l}\\
&=&{\inn{[L(\l),\nabla_{L(\l )}{\tilde F}_{m_i,i}]}{(\nabla_{L(\l )}G)_+}}_{\l}\\
                                         &=&0,
\end{eqnarray*}%
where we have used the result     $\nabla_{L(\l )}{\tilde F}_{m_i,i}\in{\tilde\g}_+$ (see item $2$ of Lemma \ref{lempft}) to
 justify the transition from second to third line
and the fact that
  ${\tilde F}_{m_i,i}$  is an  $\ad$-invariant function  on  $\tilde\g$ to obtain the last line.
\end{proof}
\begin{corollary}\label{corpf}
The  rank $\Rk(\TT_{\lambda},{\PB}_{\tilde R})$ of the  Poisson  $\tilde R$-bracket on  $\TT_{\lambda}$ is
lower or equal to  $\dim\g-\ell$.
\end{corollary}
\begin{proof}
According   to Proposition \ref{caspft}, for every $i=1,\dots,\ell$,  the  functions ${\tilde
F}_{m_i,i}$ are   Casimirs for the Poisson bracket
  ${\PB}_{\tilde R}$.
Therefore we need to show that these functions  are independent on  $\TT_{\l}$. For this,   it suffices to prove
that the   differentials with respect to the variable  $X$  of    $\tilde F_{m_i,i}$, for  $ 1\leqslant
i\leqslant \ell$ are independent. According to  the first  point of  Lemma
\ref{lempft}, for every  $1\leqslant i\leqslant \ell$ and every  $L(\l)=\l
  e_{-\b}+e+X+\l^{-1}Y$, where  $X\in\g_{\leqslant 0}$ and  $Y\in{\g_{>0}}$ , we have:
\begin{eqn*}
{\tilde F}_{m_i,i}(L(\l))=\can{\diff_Y{ P}_i,X}.
\end{eqn*}%
Then the partial derivative of  ${\tilde F}_{m_i,i}$ with respect to $X$  at
   the point $L(\l)$ is equal to
\begin{equation}\label{mlk12}
\pp{{\tilde F}_{m_i,i}}{X}(\l^{-1}Y+X+e+\l e_{-\b})=\diff_Y P_i.
\end{equation}%
In particular, at the   point  $L(\l)=\l e_{-\b}+e+X+\l^{-1}e$, (where
$e=\sum_{i=1}^{\ell}e_i$ and  $X\in\g_{\leqslant 0}$ is  arbitrary),  Equation (\ref{mlk12}) becomes:
\begin{eqn*}
\can{\pp{{\tilde F}_{m_i,i}}{X}(\l e_{-\b}+X+e+\l^{-1}e),A}=\can{\diff_e P_i,A},\qquad \forall A\in\g_{\leqslant 0}\cap T_{L(\l)}\TT_{\l}.
\end{eqn*}%
Since $e$ is  regular element of $\g$, according to  the  theorems of  Kostant
\cite[Theorem 9]{Bertram} and  \cite[Theorem  5.2]{kostant}, the
differential of the
   family  $(P_1,\dots,P_{\ell})$ are   independent at $e$.
 Moreover, since  $e\in\g_{\geq 1}$,  the  restrictions   to  $\g_{\leqslant
     0}$ of this family
are also independent because  their  gradient are in  $\g_{\geq 1}$.  Therefore the  family  $({\tilde
  F}_{m_1,1},\dots,{\tilde F}_{m_{\ell},\ell})$ is
   independent on   $\TT_{\l}$.
\end{proof}
\begin{proposition}\label{prkhim1}
The  rank  ${\Rk}({\TT_{\lambda}}, {\PB}_{\tilde R})$ of the Poisson  $\tilde R$-bracket on
$\TT_{\lambda}$ is  equal to $\dim\g-~\ell$.
\end{proposition}
According to  Corollary  \ref{corpf}, to show
Proposition  \ref{prkhim1} it suffices  to find a point
$L_0(\l)\in\TT_{\l}$ where  the rank of the Poisson structure  is  $\dim\g-\ell$.
We start by stating  Lemma \ref{kh77}, the proof of which is a direct computation describing explicitly the Poisson structure of $\TT_\l$. Notice that, although $\TT_\l$ is an affine subspace of $\tilde \g$, the Poisson structure obtained by restriction to $\TT_\l$ is linear.
\begin{lemma}\label{kh77}
   For all $i=1,\dots,\ell$ and all $ \a\in\Phi_+$, let
 $x_i,x_{-\alpha},y_{\alpha}$  be the coordinates functions  on  $\TT_{\l}$,   defined at  every point
 $L(\l)=\l e_{-\b}+e+X+\l^{-1}Y$ of   $\TT_{\l}$,  where  $X\in\g_{\leqslant 0}$ and  $Y\in\g_{>0}$,  by:
\begin{eqn*}
\left\{
\renewcommand{\arraystretch}{0.9}
\begin{array}{llll}
\can{x_i,L(\l)}&:=&\inn{h_i}{X},&\\
\can{x_{-\alpha},L(\l)}&:=&\inn{e_{\alpha}}{X}, &\\
\can{y_{\alpha},L(\l)}&:=&\inn{ e_{-\alpha}}{Y},&
\end{array}
\right.
\end{eqn*}%
 The expression of  the Poisson  $\tilde R$-bracket on  $\TT_{\lambda}$ is given, for every
 $ 1\leqslant i,j\leqslant \ell$  and every  $\alpha,\gamma\in\Phi_+$, by:
\begin{equation}\label{poisson_crochet}
\left\{
\begin{array}{ll}
\renewcommand{\arraystretch}{0.9}
{\pb{x_i,x_j}}_{\tilde R}=0,& \\
{\pb{x_i,x_{-\alpha}}}_{\tilde R}=\alpha(h_i)x_{-\alpha},&\\
{\pb{x_i,y_{\alpha}}}_{\tilde R}=-\alpha(h_i)y_{\alpha},&\\
{\pb{x_{-\alpha},x_{-\gamma}}}_{\tilde  R}=\eta_{\alpha+\gamma}N_{\alpha,\gamma}x_{-\alpha-\gamma},&
\\
{\pb{x_{-\alpha},y_{\gamma}}}_{\tilde R}=\eta_{\gamma-\a}N_{\alpha,-\gamma}y_{\gamma-\alpha},&\\
{\pb{y_{\alpha},y_{\gamma}}}_{\tilde R}=0,&
\end{array}
\right.
\end{equation}%
where
$
\eta_{\a}=
\left\{
\renewcommand{\arraystretch}{0.9}
\begin{array}{ll}
1, &\textrm{ if }\a\in\Phi_+, \\
0,&\textrm{ otherwise,}
\end{array}
\right.
$
and   $N_{\a\gamma}=\pm(p+1)$, with  $p:=max \{ n \mid \gamma-n\a\in\Phi\}$.
\end{lemma}
We now show  Proposition \ref{prkhim1}.
 \begin{proof}
 Let    $b_1,\dots,b_{\ell}$ be  non-zero  constants and let
 \begin{equation}\label{eqL0}
L_0(\lambda):=\sum_{i=1}^{\ell}(1+\lambda^{-1} b_i)e_i +\lambda e_{-\beta}.
\end{equation} %
 According   to (\ref{poisson_crochet}), for every  $ 1\leqslant i,j\leqslant \ell $,  the  Poisson  $\tilde R$-bracket at the
 point $L_0(\lambda)$  is given by:
\begin{equation}\label{vchss1}
\left\{
\renewcommand{\arraystretch}{0.9}
\begin{array}{ll}
{\pb{x_i,x_j}}_{\tilde R}=0,&  \\
{\pb{x_i,x_{-\alpha}}}_{\tilde R}=0,&\\
{\pb{x_i,y_{\alpha}}}_{\tilde R}=\left\{
                                  \begin{array}{ll}
                                   -c_{ji}b_j&~~ \textrm{if }
                                   \alpha\textrm{ is a simple root
                                     }\alpha _j,\\
                                   0&~~\textrm{otherwise,}
                                  \end{array}
                                  \right.&\\
{\pb{x_{-\alpha},x_{-\gamma}}}_{\tilde R}=0,&\\
{\pb{x_{-\alpha},y_{\gamma}}}_{\tilde R}=\left\{
                                       \begin{array}{ll}
                                       N_{\alpha,-\gamma}b_i
                                       &~~\textrm{if  }
                                       \gamma-\a\textrm{ is a simple root
                                          }
                               \a_i, \\
                                       0 &~~\textrm{otherwise,}
                                     \end{array}
                                      \right. &\\
{\pb{y_{\alpha},y_{\gamma}}}_{\tilde R}=0,&
\end{array}
\right.
\end{equation}%
where $(c_{ij})_{ 1\leqslant i,j\leqslant \ell}$ is the Cartan
matrix of $\g$. We denote by
$\gamma_1,\dots,\gamma_{\frac{\dim\g-\ell}{2}}$ the   positive roots
of  $\g$ and we choose the indices such that    $|\gamma_1|\leqslant |\gamma_2|\leqslant \dots\leqslant
|\gamma_{\frac{\dim\g-\ell}{2}}|$. It will be convenient to denote by 
  $(z_1,\dots,z_{\dim\g})$  the system  of  coordinates  given by:
\begin{eqn*}
\renewcommand{\arraystretch}{1.8}
\left\{
\renewcommand{\arraystretch}{0.9}
\begin{array}{lll}
z_i=x_i,&& 1\leqslant i\leqslant \ell,\\
z_{\ell+k}=x_{-\gamma_k},& & 1\leqslant k\leqslant\frac{\dim \g-\ell}{2},\\
z_{(\frac{\dim\g+\ell}{2}+j)}=y_{\gamma_j},&& 1\leqslant j\leqslant \frac{\dim\g-\ell}{2}.
\end{array}
\right.
\end{eqn*}%
 By using the formulas of system (\ref{vchss1}), one  establishes   the matrix
 $M=({\pb{z_i,z_j}}_{\tilde R})_{ 1\leqslant i,j\leqslant \dim\g}$  of  the Poisson  $\tilde
R$-bracket  computed  at  the  point
$L_0(\lambda)$ given in  (\ref{eqL0}). We obtain a matrix of the form
\begin{equation}
M=\left(
\begin{array}{cc}
0&-{\Lambda}^{T}\\
\Lambda &0
\end{array}
\right),
\end{equation}%
where  $\Lambda$ is the following  block diagonal matrix of size
$\frac12(\dim\g-\ell)\times\frac12(\dim\g+\ell)$
\begin{equation}
\Lambda=\left(
\renewcommand{\arraystretch}{0.9}
\begin{array}{ccccc}
\Lambda_0&         &         &0              &0 \\
         &\Lambda_1&         &               &\vdots \\
         &         &\ddots   &               &\vdots\\
  0      &         &         &\Lambda_{{m_{\ell}-1}}&0
\end{array}
\right),
\end{equation}%
where the $0$ aligned vertically
 at  right hind  of the matrix represents a single column and not a group of columns, and  $\Lambda_0\dots,\Lambda_{{m_{\ell}-1}}$
are matrices whose expressions shall be given later.\\
 Let  $B=\left(\begin{array}{ccc}
                                 b_1& &0\\
                                 &\ddots&\\
                                0&&  b_{\ell}
                              \end{array}
                           \right)$ and  $C=(c_{ij})_{1\leqslant i,j\leqslant \ell}$  be the Cartan  matrix of
                           $\g$, we have   $\Lambda_0=BC$.\\
  We recall that
\begin{eqn*}
\left\{
\renewcommand{\arraystretch}{0.5}
\begin{array}{l}
\dim\g_0=\dim{\g_1}=\dim\g_{-1}=\ell,\\
\dim\g_{m_{\ell}}=1,\\
\sum_{i=1}^{m_{\ell}}\dim\g_i=\frac{1}{2}(\dim\g-\ell).
\end{array}
\right.
\end{eqn*}%
 We denote by $d_i$ the dimension of $\g_i$ and we denote, for  $k\neq 0$, by   $(\gamma_1,\dots,\gamma_{d_k})$
 a  basis of roots of
$\g$ of length $k$, $(\beta_1,\dots,\beta_{d_{k+1}})$ a  basis  of
 roots of  $\g$   of length $k+1$. By definition:
\begin{equation}\label{matrgamma}
\Lambda _k
=\left(
\renewcommand{\arraystretch}{0.9}
           \begin{array}{ccc}
\X_{y_{\beta_1}}[x_{-\gamma_1}] &  \cdots &\X_{y_{\beta_{d_{k+1}}}}[x_{-\gamma_1}]\\
    \vdots                 &         &\vdots\\
\X_{y_{\beta_1}}[x_{-\gamma_{d_{k}}}]&\cdots& \X_{y_{\beta_{d_{k+1}}}}[x_{-\gamma_{d_{k}}}]
          \end{array}
           \right)^T.
\end{equation}%
 To show 
$\Rk(L_0(\lambda),{\PB}_{\tilde R})=\dim\g-\ell$ it is necessary and sufficient to prove that the rank of
matrix $\Lambda$ is  $\frac12(\dim\g-\ell)$. In turn this   is equivalent to
show that first the rank of  $\Lambda _0$, is  $\ell$ and  that  every matrix $ \Lambda_k $, for $1\leqslant k\leqslant m_{\ell}-1
$ is of rank $d_{k+1}$.
\\
\emph{(1)} The Cartan matrix is invertible, and assuming that
  $b_1,\dots, b_{\ell}$ are   non-zero,  the  matrix $\Lambda_0=BC$ is
  invertible also so that the  rank of $\Lambda_0$ is  $\ell$.\\
\emph{(2)}
 We recall  that, for every  $1\leqslant i \leqslant d_k$ and for every
  $1\leqslant j \leqslant d_{k+1}$, we have:
\begin{eqn*}
\X_{y_{{\b_{j}}}}[x_{-\gamma_i}]=
\left\{
\begin{array}{ll}
N_{-\b_j,\gamma_i}b_{p},& \textrm{ if  } \b_j-\gamma_i \textrm{ is a simple
  root } \a_p,\\
0,&\textrm{ otherwise.}
\end{array}
\right.
\end{eqn*}%
Let  $1\leqslant j\leqslant d_{k+1}$.
 For every  $\b_j$, there  exists a
    index $i\in\{1,\dots, d_k\}$ and a
    index  $F(i,j)\in\{1, \dots,\ell\}$, such that:
\begin{eqn*}
\b_j=\gamma_i+\a_{F(i,j)}.
\end{eqn*}%
This  implies that:
\begin{eqn*}
 \X_{y_{\b_j}}[x_{\gamma_i}]=N_{-\b_j,\gamma_i}b_{F(i,j)}.
\end{eqn*}%
By  construction, the  above constant $N_{-\b_j,\gamma_i}$ is   non-zero  and  equal to   $1$.
We prove, for each simple Lie algebra, for
$b_1,\dots,b_{\ell}$   generic, the rank of the matrix
 $\Lambda_k$ is  $d_{k+1}$, for every  $k=1,\dots,m_{\ell}-1$.\\
\emph{(a)} To prove  the result for the classical simple Lie algebras of  $\g$  of  type $A_{\ell}, B_{\ell},C_{\ell}$ 
and  $D_{\ell}$,  we fix an order on the roots of the same length. 
Then we show that the matrices henceforth  obtained have the required  rank.

{\bf{ Case $A_{\ell}$: }}
 Let  $\g$ be the simple Lie algebra  of type  $A_{\ell}$  and  let
  $\a_1,\dots\a_{\ell}$ be the simple roots of  $\g$. We choose to
 arrange   the roots of
  length $k$ of $\g$ in the following  (lexicographic) order $
\gamma_1=\a_1+\dots+\a_k,
\gamma_2=\a_2+\dots+\a_{k+1},\dots,
\gamma_{{\ell}-k}=\a_{{\ell}-k}+\dots+\a_{{\ell}-1},
\gamma_{{\ell}-k+1}=\a_{\ell-k+1}+\dots+\a_{\ell},
$
and the roots of $\g$ of length $k+1$ in lexicographic order, which gives the  array below where all the decompositions
of a root of length $k+1$ as a sum of a simple root with a root of length $k$ and we have, for every  $1\leqslant j\leqslant \ell-k$,
\begin{equation}
\b_{j}=\gamma_{j}+\a_{k+j}=\gamma_{j+1}+\a_j.
\end{equation}%
 The  matrix  $\Lambda_k^T$, defined in  (\ref{matrgamma}) is   of  the form:
\begin{equation}
\Lambda_k^T=\left(
\renewcommand{\arraystretch}{0.9}
\begin{array}{ccccc}
b_{k+1}&       &       &      & 0 \\
b_1   &b_{k+2} &       &      &   \\
      &b_2     &\ddots &      &   \\
      &        &\ddots&\ddots &   \\
      &        &      &\ddots &b_{\ell}  \\
0     &        &      &       &b_{\ell-k}
\end{array}
\right).
\end{equation}%
By removing the last line of  $\Lambda_k^T$, we obtain a lower triangular square  ($d_{k+1}\times d_{k+1}$)
matrix   $\Gamma_k$,  which is of  rank $d_{k+1}$, when $b_{k+1},\dots b_{\ell}$ are all non-zero. This implies
that the rank of $\Lambda_k $ is $d_{k +1}$.

{\bf{Cass $B_{\ell}$: }}
Let  $\g$ be a simple Lie algebra of  type $B_{\ell}$ and let $(\a_1,\dots,\a_{\ell})$ a basis of simple roots of  $\g$. The positive roots of $\g$ have the following expressions
\begin{eqn*}
\left\{
\renewcommand{\arraystretch}{0.9}
\begin{array}{ll}
\l_i=\a_i+\dots+\a_{\ell},&\qquad 1\leqslant i\leqslant\ell,\\
\l_i-\l_j=\a_i+\dots+\a_{j-1},&\qquad 1\leqslant i<j\leqslant \ell,\\
\l_i+\l_j=\a_i+\dots +\a_{j-1}+2(\a_j+\dots+\a_{\ell}),&\qquad 1\leqslant
i<j\leqslant \ell.
\end{array}
\right.
\end{eqn*}%
To establish the rank of the matrix $\Lambda_k$, we need to discuss following the parity of $k$. For $k$ even, we choose to arrange the roots of $\g$ of length $k$  in lexicographic order (lexicographic with respect to  $(\l_1,\dots,\l_\ell)$), to wit
$
\gamma_1=\l_1-\l_{k+1},\dots,\gamma_{\ell-k}=\l_{\ell-k}-\l_{\ell}, \gamma_{\ell-k+1}=\l_{\ell-k+1}, \gamma_{\ell-k+2}=\l_{\ell-k+2}+\l_{\ell},\dots, \gamma_{\ell-\frac{k}{2}}=\l_{\ell-\frac{k}{2}-1}+\l_{\ell-\frac{k}{2}+3},
\gamma_{\ell-\frac{k}{2}}=\l_{\ell-\frac{k}{2}}+\l_{\ell-\frac{k}{2}+2}
$
and the roots of $\g$ of length $k+1$ in lexicographic order, which gives the  array below where all the decompositions
of a root of length $k+1$ as a sum of a simple root with a root of length $k$ have been indicated on the right column:
\begin{eqn*}
\renewcommand{\arraystretch}{0.9}
\begin{array}{lllll}
\b_1&=&\l_1-\l_{k+2}&=&\left\{\begin{array}{l}\a_1+\gamma_2,\\
                                       \gamma_1+\a_{k+1},
                                        \end{array}
                                          \right.\\
\vdots & &\vdots  & & \\
\b_{\ell-k-1}&=&\l_{\ell-k-1}-\l_{\ell}&=&\left\{ \begin{array}{l} \a_{\ell-k-1}+\gamma_{l-k},\\
                                                             \gamma_{\ell-k-1}+\a_{\ell-1},
                                              \end{array}
                                     \right.\\
\b_{\ell-k}&=&\l_{\ell-k}&=&\left\{ \begin{array}{l} \a_{\ell-k}+\l_{\ell-k+1},\\
                                                 \gamma_{\ell-k}+\a_{\ell},
                                \end{array} \right.\\
\b_{\ell-k+1}&=&\l_{\ell-k+1}+\l_{\ell}&=&\left\{\begin{array}{l} \gamma_{\ell-k+1}+\a_{\ell},\\
                                                            \a_{\ell-k+1}+\gamma_{\ell-k+2},
                                            \end{array}\right.\\
\vdots&&\vdots &&\\
\b_{\ell-\frac{k}{2}-1}&=&\l_{\ell-\frac{k}{2}-1}+\l_{\ell-\frac{k}{2}+2}
                                                                  &=&\left\{\begin{array}{l}
                                                                 \a_{\ell-\frac{k}{2}-1}+\gamma_{\ell-\frac{k}{2}},\\
                                                                \gamma_{\ell-\frac{k}{2}-1}+\a_{\ell-\frac{k}{2}+2},
                                                                 \end{array}
                                                                 \right.\\
\b_{\ell-\frac{k}{2}}&=&\l_{\ell-\frac{k}{2}}+\l_{\ell-\frac{k}{2}+1}&=&\gamma_{\ell-\frac{k}{2}}+\a_{\ell-\frac{k}{2}+1}.
\end{array}
\end{eqn*}%
In view of the previous array, the  matrix $\Lambda_k^T$, defined in (\ref{matrgamma}) takes the following form:
\begin{eqn*}
\Lambda_k^T=\left(
\renewcommand{\arraystretch}{0.9}
\begin{array}{cccccccccc}
b_{k+1} &      &         &          &            &            &      &&&\\
b_1    &\ddots &         &          &           &            &   0     &&&\\
       &\ddots &\ddots   &          &           &            &        &&&\\
       &       &b_{\ell-k-1} &b_{\ell}  &           &            &        &&&\\
       &       &         &b_{\ell-k}&b_{\ell}   &            &        &&&\\
       &       &         &         & b_{\ell-k+1}& b_{\ell-1} &        &&&\\
       &       &         &         &            &b_{\ell-k+2} &\ddots &&&\\
       &       &         &         &            &            &\ddots &\ddots &&\\
      &0       &         &         &            &            &       &\ddots&b_{\ell-\frac{k}{2}+2}&\\
       &       &         &         &            &            &       &&b_{l-\frac{k}{2}-1}&b_{\ell-\frac{k}{2}+1}
\end{array}
\right).
\end{eqn*}%
We notice that  $\Lambda_k^T$ is a  lower  triangular square $(d_{k+1}\times d_{k+1})$ matrix.
 Its determinant is a product of a finite number of
$b_i$, therefore it is non-zero (we recall that the
$b_1,\dots,b_{\ell}$ all different from zero). This implies that the rank of
$\Lambda_k^T$ is $d_{k+1}$.

For  $k$ odd, we arrange the roots of $\g$ of lengths  $k$ in lexicographic order, to wit
$\gamma_1=\l_1-\l_{k+1},\dots,
\gamma_{\ell-k}=\l_{\ell-k}-\l_{\ell},
\gamma_{\ell-k+1}=\l_{\ell-k+1},
\gamma_{\ell-k+2}=\l_{\ell-k+2}+\l_{\ell},\dots,
\gamma_{\ell-\frac{k-1}{2}-1}=
\l_{\ell-\frac{k-1}{2}-1}+\l_{\ell-\frac{k-1}{2}+2},
\gamma_{\ell-\frac{k-1}{2}}=\l_{\ell-\frac{k-1}{2}}+\l_{\ell-\frac{k-1}{2}+1}
$
and the roots of $\g$ of length $k+1$ in lexicographic order, which gives the  array below where all the decompositions
of a root of length $k+1$ as a sum of a simple root with a root of length $k$ have been indicated on the right column:
\begin{eqn*}
\renewcommand{\arraystretch}{0.9}
\begin{array}{lllll}
\b_1&=&\l_1-\l_{k+2}&=&\left\{\begin{array}{l}\a_1+\gamma_2,\\
                                       \gamma_1+\a_{k+1},
                                        \end{array}
                                          \right.\\
\vdots & &\vdots   && \\
\b_{\ell-k-1}&=&\l_{\ell-k-1}-\l_{\ell}&=&\left\{ \begin{array}{l} \a_{\ell-k-1}+\gamma_{l-k},\\
                                                             \gamma_{\ell-k-1}+\a_{\ell-1},
                                              \end{array}
                                     \right.\\
\b_{\ell-k}&=&\l_{\ell-k}&=&\left\{ \begin{array}{l} \a_{\ell-k}+\l_{\ell-k+1},\\
                                                 \gamma_{\ell-k}+\a_{\ell},
                                \end{array} \right.\\
\b_{\ell-k+1}&=&\l_{\ell-k+1}+\l_{\ell}&=&\left\{\begin{array}{l} \gamma_{\ell-k+1}+\a_{\ell},\\
                                                            \a_{\ell-k+1}+\gamma_{\ell-k+2},
                                            \end{array}\right.\\
\vdots&&\vdots &&\\
\b_{\ell-\frac{k-1}{2}-1}&=&\l_{\ell-\frac{k-1}{2}-1}+\l_{\ell-\frac{k-1}{2}+1}
                                                                  &=&\left\{\begin{array}{l}
                                                                 \a_{\ell-\frac{k-1}{2}-1}+\gamma_{\ell-\frac{k-1}{2}},\\
                                                                \gamma_{\ell-\frac{k-1}{2}-1}+\a_{\ell-\frac{k-1}{2}+1}.
                                                                 \end{array}
                                                                 \right.
\end{array}
\end{eqn*}%
In view of the previous array, the  matrix $\Lambda_k^T$, defined in (\ref{matrgamma}) takes the following form:
have the following  form:
\begin{eqn*}
\Lambda_k^T=\left(
\renewcommand{\arraystretch}{0.9}
\begin{array}{ccccccccc}
b_{k+1}&      &        &       &          &         &       0     &&\\
b_1    &\ddots &      &       &           &         &           &&\\
       &\ddots&\ddots &       &           &        &            &&\\
       &      &\ddots &b_{\ell-1}&         &        &            &&\\
       &      &       &b_{\ell-k-1}&b_{\ell}&        &            &&\\
       &      &       &          &b_{\ell-k}&b_{\ell}&            &&\\
       &      &       &          &         &b_{\ell-k+1}&b_{\ell-1}&&\\
       &      &       &          &         &           &\ddots   &\ddots&\\
       &   0   &       &          &         &           &         &\ddots &b_{\ell-\frac{k-1}{2}+1}\\
       &      &       &          &         &          &          &   &b_{\ell-\frac{k-1}{2}-1}
\end{array}
\right).
\end{eqn*}%
By removing the last line of  $\Lambda_k^T$, defined in (\ref{matrgamma}), we  obtain a lower triangular square  ($d_{k+1}\times d_{k+1}$)
matrix   $\Gamma_k$ and which is of  rank $d_{k+1}$, when $b_{k+1},\dots b_{\ell}$ are  all non-zero. This implies
that the rank of $\Lambda_k $ is $d_{k +1}$.

{\bf{Case $C_{\ell}$: }}
Let  $\g$ be a simple Lie algebra of   type $C_{\ell}$ and  let
$(\a_1,\dots,\a_{\ell})$ be a basis of simple roots of   $\g$. The expressions of the positive roots of
 $\g$ are
\begin{eqn*}
\left\{
\renewcommand{\arraystretch}{0.9}
\begin{array}{ll}
2\l_i=2(\a_i+\dots+\a_{\ell-1})+\a_{\ell},&\qquad 1\leqslant i\leqslant \ell,\\
\l_i-\l_j=\a_i+\dots+\a_{j-1},&\qquad 1\leqslant i<j\leqslant \ell,\\
\l_i+\l_j=\a_i+\dots+\a_{j-1}+2(\a_j+\dots+\a_{\ell-1})+\a_{\ell},&\qquad
1\leqslant i<j\leqslant\ell.
\end{array}
\right.
\end{eqn*}%
To compute the rank of the matrix
$\Lambda_k$, we discuss following the parity of $k$. For $ k $ even,
we choose to arrange the roots of $\g$ of length $k$ in lexicographic order, to wit
$
\gamma_1=\l_1-\l_{k+1}, \dots,
\gamma_{\ell-k-1}=\l_{\ell-k-1}-\l_{\ell-1},
\gamma_{\ell-k}=\l_{\ell-k}-\l_{\ell},
\gamma_{\ell-k+1}=\l_{\ell-k+1}+\l_{\ell},
\gamma_{\ell-k+2}=\l_{\ell-k+2}+\l_{\ell-1},\dots,
\gamma_{\ell-\frac{k}{2}-1}=\l_{\ell-\frac{k}{2}-1}+\l_{\ell-\frac{k}{2}+2},
\gamma_{\ell-\frac{k}{2}}=\l_{\ell-\frac{k}{2}}+\l_{\ell-\frac{k}{2}+1},
$
and the roots of $\g$ of length $k+1$ in lexicographic order, which gives the  array below where all the decompositions
of a root of length $k+1$ as a sum of a simple root with a root of length $k$ have been indicated on the right column:
\begin{eqn*}
\renewcommand{\arraystretch}{0.9}
\begin{array}{lllll}
\b_1&=&\l_1-\l_{k+2}&=&\left\{\begin{array}{l}\a_1+\gamma_2,\\
                                       \gamma_1+\a_{k+1},
                                        \end{array}
                                          \right.\\
\vdots && \vdots && \\
\b_{\ell-k-1}&=&\l_{\ell-k-1}-\l_{\ell}&=&\left\{\begin{array}{l}\a_{\ell-k-1}+\gamma_{\ell-k},\\
                                       \gamma_{\ell-k-1}+\a_{\ell-1},
                                        \end{array}
                                          \right.\\
\b_{\ell-k}&=&\l_{\ell-k}+\l_{\ell}&=&\left\{\begin{array}{l}\a_{\ell}+\gamma_{\ell-k},\\
                                       \gamma_{\ell-k+1}+\a_{\ell-k},
                                        \end{array}
                                          \right.\\
\b_{\ell-k+1}&=&\l_{\ell-k+1}+\l_{\ell-1}&=&\left\{\begin{array}{l}\a_{\ell-1}+\gamma_{\ell-k+1},\\
                                       \gamma_{\ell-k+2}+\a_{\ell-k+1},
                                        \end{array}
                                          \right.\\
\vdots  &&\vdots && \\
\b_{\ell-\frac{k}{2}-1}&=&\l_{\ell-\frac{k}{2}-1}+\l_{\ell-\frac{k}{2}+1}&=&\left\{\begin{array}{l}
                                       \a_{\ell-\frac{k}{2}+1}+\gamma_{\ell-\frac{k}{2}-1},\\
                                       \gamma_{\ell-\frac{k}{2}}+\a_{\ell-\frac{k}{2}-1},
                                        \end{array}
                                          \right.\\
\b_{\ell-\frac{k}{2}}&=&2\l_{\ell-\frac{k}{2}}&=&\a_{\ell-\frac{k}{2}}+\gamma_{\ell-\frac{k}{2}}.
\end{array}%
\end{eqn*}%
Therefore the matrix  $\Lambda_k^T$, defined in  (\ref{matrgamma}) has the following form:
\begin{eqn*}
\Lambda_k^T=\left(
\renewcommand{\arraystretch}{0.9}
\begin{array}{cccccccccc}
b_{k+1}& & & & & & & & & \\
b_1   &\ddots& & & & & &0 & &\\
      &\ddots&\ddots& & & & & & &\\
      &      &\ddots&b_{\ell-1}& & & & & & \\
      &      &      &b_{\ell-k-1}&b_{\ell}& & & & & \\
      &      &      &           &b_{\ell-k}&b_{\ell-1} &  &   & &\\
      &      &      &           &         &b_{\ell-k+1}&\ddots& & & \\
      &      &      &           &         &           &\ddots&\ddots &&\\
     0 &      &      &           &         &           &      &\ddots&b_{\ell-\frac{k}{2}+1}&\\
      &      &      &           &         &           &      &      &b_{\ell-\frac{k}{2}-1}&b_{\ell-\frac{k}{2}}
\end{array}
\right).
\end{eqn*}%
We notice that  $\Lambda_k^T$ is a lower triangular square $(d_{k+1}\times d_{k+1})$   matrix
$(d_{k+1}\times d_{k+1})$.   Its determinant is a product of a finite number of
$b_i$, therefore  it is non-zero. This implies that the rank of $\Lambda_k^T$ is $d_{k+1}$.

We consider now  the case where $k$ is odd. The roots of $\g$ of length
$k$ are ordered by lexicographic order, to wit
$
\gamma_1=\l_1-\l_{k+1}, \dots,
\gamma_{\ell-k-1}=\l_{\ell-k-1}-\l_{\ell-1},
\gamma_{\ell-k}=\l_{\ell-k}-\l_{\ell},
\gamma_{\ell-k+1}=\l_{\ell-k+1}+\l_{\ell},
\gamma_{\ell-k+2}=\l_{\ell-k+2}+\l_{\ell-1},
\dots,
\gamma_{\ell-\frac{k-1}{2}-2}=\l_{\ell-\frac{k-1}{2}-2}+\l_{\ell-\frac{k-1}{2}+2},
\gamma_{\ell-\frac{k-1}{2}-1}=\l_{\ell-\frac{k-1}{2}-1}+\l_{\ell-\frac{k-1}{2}+1},
\gamma_{\ell-\frac{k-1}{2}}=2\l_{\ell-\frac{k-1}{2}},
$
and the roots of $\g$ of length $k+1$ in lexicographic order, which gives the  array below where all the decompositions
of a root of length $k+1$ as a sum of a simple root with a root of length $k$ have been indicated on the right column:
\begin{eqn*}
\renewcommand{\arraystretch}{0.9}
\begin{array}{lllll}
\b_1&=&\l_1-\l_{k+2}&=&\left\{\begin{array}{l}\a_1+\gamma_2,\\
                                       \gamma_1+\a_{k+1},
                                        \end{array}
                                          \right.\\
\vdots && \vdots \\
\b_{\ell-k-1}&=&\l_{\ell-k-1}-\l_{\ell}&=&\left\{\begin{array}{l}\a_{\ell-k-1}+\gamma_{\ell-k},\\
                                       \gamma_{\ell-k-1}+\a_{\ell-1},
                                        \end{array}
                                          \right.\\
\b_{\ell-k}&=&\l_{\ell-k}+\l_{\ell}&=&\left\{\begin{array}{l}\a_{\ell}+\gamma_{\ell-k},\\
                                       \gamma_{\ell-k+1}+\a_{\ell-k},
                                        \end{array}
                                          \right.\\
\b_{\ell-k+1}&=&\l_{\ell-k+1}+\l_{\ell-1}&=&\left\{\begin{array}{l}\a_{\ell-1}+\gamma_{\ell-k+1},\\
                                       \gamma_{\ell-k+2}+\a_{\ell-k+1},
                                        \end{array}
                                          \right.\\
\vdots  &&\vdots\\
\b_{\ell-\frac{k-1}{2}-2}&=&\l_{\ell-\frac{k-1}{2}-2}+\l_{\ell-\frac{k-1}{2}+1}&=&\left\{\begin{array}{l}
                                       \a_{\ell-\frac{k-1}{2}+1}+\gamma_{\ell-\frac{k-1}{2}-2},\\
                                       \gamma_{\ell-\frac{k-1}{2}-1}+\a_{\ell-\frac{k-1}{2}-2},
                                        \end{array}
                                          \right.\\
\b_{\ell-\frac{k-1}{2}-1}&=&\l_{\ell-\frac{k-1}{2}-1}+\l_{\ell-\frac{k-1}{2}}&=&\left\{\begin{array}{l}
                                       \a_{\ell-\frac{k-1}{2}}+\gamma_{\ell-\frac{k-1}{2}-1},\\
                                       \gamma_{\ell-\frac{k-1}{2}}+\a_{\ell-\frac{k-1}{2}-1}.
                                        \end{array}
                                          \right.
\end{array}%
\end{eqn*}%
Therefore the matrix  $\Lambda_k^T$ defined in (\ref{matrgamma}) takes the following form:
\begin{eqn*}
\Lambda_k^T=\left(
\renewcommand{\arraystretch}{0.9}
\begin{array}{cccccccccc}
b_{k+1}& & & & & & & &0 & \\
b_1   &\ddots& & & & & & & &\\
      &\ddots&\ddots& & & & & & &\\
      &      &\ddots&b_{\ell-1}& & & & & & \\
      &      &      &b_{\ell-k-1}&b_{\ell}& & & & & \\
      &      &      &           &b_{\ell-k}&b_{\ell-1} &  &   & &\\
      &      &      &           &         &b_{\ell-k+1}&\ddots& & & \\
      &      &      &           &         &           &\ddots&\ddots &&\\
      &      &      &           &         &           &      &\ddots&b_{\ell-\frac{k-1}{2}+1}&\\
    0  &      &      &           &         &           &      &      &b_{\ell-\frac{k-1}{2}-2}&b_{\ell-\frac{k-1}{2}}\\
      &      &      &           &         &           &      &      &&b_{\ell-\frac{k-1}{2}-1}
\end{array}
\right).
\end{eqn*}%
By removing the last line of  $\Lambda_k^T$, we  obtain a lower triangular square  ($d_{k+1}\times d_{k+1}$)
matrix   $\Gamma_k$ which is of  rank $d_{k+1}$ when $b_{k+1},\dots b_{\ell}$ are all non-zero. This implies
that the rank of $\Lambda_k $ is $d_{k +1}$.

{\bf{Case $D_{\ell}$: }}
Let  $\g$ be a simple Lie algebra  of type $D_{\ell}$ and let $(\a_1,\dots,\a_{\ell})$ be a basis of simple roots of $\g$. The positive roots of $\g$ are
\begin{eqn*}
\left\{
\renewcommand{\arraystretch}{0.9}
\begin{array}{ll}
\l_i-\l_j=\a_i+\dots+\a_{j-1},&\qquad 1\leqslant i<j\leqslant \ell,\\
\l_i+\l_{\ell}=\a_i+\dots+\a_{\ell-2}+\a_{\ell},&\qquad 1\leqslant
i<\ell,\\
\l_i+\l_j=\a_i+\dots+\a_{j-1}+2(\a_{j}+\dots+\a_{\ell-2})+\a_{\ell-1}+\a_{\ell},&\qquad
1\leqslant i<j<\ell.
\end{array}
\right.
\end{eqn*}%
As in the case of $B_{\ell}$, to calculate the rank of the matrix $\Lambda_k$ we study separately the cases where the integer $k$ is even and odd. Let us start with the case $ k $ is even. We arrange the roots of $\g$ of length $k$ in lexicographic order, to wit:
$
\gamma_1=\l_1-\l_{k+1},\dots,
\gamma_{\ell-k-1}=\l_{\ell-k-1}-\l_{\ell-1},
\gamma_{\ell-k}=\l_{\ell-k}-\l_{\ell},
\gamma_{\ell-k+1}=\l_{\ell-k}+\l_{\ell},
\gamma_{\ell-k+2}=\l_{\ell-k+1}+\l_{\ell-1},\dots,
\gamma_{\ell-\frac{k}{2}-1}=\l_{\ell-\frac{k}{2}-2}+\l_{\ell-\frac{k}{2}+2},
\gamma_{\ell-\frac{k}{2}}=\l_{\ell-\frac{k}{2}-1}-\l_{\ell-\frac{k}{2}+1},
$
and the roots of $\g$ of length $k+1$ in lexicographic order, which gives the  array below where all the decompositions
of a root of length $k+1$ as a sum of a simple root with a root of length $k$ have been indicated on the right column:
\begin{eqn*}
\renewcommand{\arraystretch}{0.9}
\begin{array}{lllll}
\b_1&=&\l_1-\l_{k+2}&=&\left\{\begin{array}{l}
                                       \a_1+\gamma_2,\\
                                       \gamma_1+\a_{k+1},
                                        \end{array}
                                          \right.\\
\vdots&&\vdots&&\\
\b_{\ell-k-1}&=&\l_{\ell-k-1}-\l_{\ell}&=&\left\{\begin{array}{l}
                                       \a_{\ell-k-1}+\gamma_{\ell-k},\\
                                       \gamma_{\ell-k-1}+\a_{\ell-1},
                                        \end{array}
                                          \right.\\
\b_{\ell-k}&=&\l_{\ell-k-1}+\l_{\ell}&=&\left\{\begin{array}{l}
                                       \a_{\ell-k-1}+\gamma_{\ell-k+1},\\
                                       \gamma_{\ell-k-1}+\a_{\ell},
                                        \end{array}
                                          \right.\\
\b_{\ell-k+1}&=&\l_{\ell-k}+\l_{\ell-1}&=&\left\{\begin{array}{l}
                                       \a_{\ell}+\gamma_{\ell-k},\\
                                       \gamma_{\ell-k+1}+\a_{\ell-1},\\
                                        \gamma_{\ell-k+2}+\a_{\ell-k},
                                        \end{array}
                                          \right.\\
\b_{\ell-k+2}&=&\l_{\ell-k+1}+\l_{\ell-2}&=&\left\{\begin{array}{l}
                                       \a_{\ell-k+1}+\gamma_{\ell-k+3},\\
                                        \gamma_{\ell-k+2}+\a_{\ell-2},
                                        \end{array}
                                          \right.\\
\vdots      &&\vdots&& \\
\b_{\ell-\frac{k}{2}-1}&=&\l_{\ell-\frac{k}{2}-2}+\l_{\ell-\frac{k}{2}+1}&=&\left\{\begin{array}{l}
                                       \a_{\ell-\frac{k}{2}-2}+\gamma_{\ell-\frac{k}{2}},\\
                                        \gamma_{\ell-\frac{k}{2}-1}+\a_{\ell-\frac{k}{2}+1},
                                        \end{array}
                                          \right.\\
\b_{\ell-\frac{k}{2}}&=&\l_{\ell-\frac{k}{2}-1}+\l_{\ell-\frac{k}{2}}&=&\a_{\ell-\frac{k}{2}}+\gamma_{\ell-\frac{k}{2}}.
\end{array}%
\end{eqn*}%
Then the matrix  $\Lambda_k^T$ defined in (\ref{matrgamma})
takes the  following form:
\begin{eqn*}
\Lambda_k^T=\left(
\renewcommand{\arraystretch}{0.9}
\begin{array}{cccccccccc}
b_{k+1}&      &      & & & & & & & \\
b_1   &\ddots& & & & & & & & \\
      &\ddots&\ddots& & & & & & & \\
      &      &\ddots&b_{\ell-1}  &b_{\ell}& & & & & \\
      &      &      &b_{\ell-k-1}&0          &b_{\ell}& & & & \\
      &      &      &           &b_{\ell-k-1}&b_{\ell-1}& & & & \\
      &      &      &           &           &b_{\ell-k}&b_{\ell-2}& & &\\
      &      &      &           &           &         &b_{\ell-k+1}&\ddots&&\\
      &      &      &           &           &         &           &\ddots&b_{\ell-\frac{k}{2}+1}&\\
      &      &      &           &           &         &           &      &b_{\ell-\frac{k}{2}-2}&b_{\ell-\frac{k}{2}}
\end{array}
\right).
\end{eqn*}%
The  matrix  $\Lambda_k^T$ is a square matrix and we verify that
\begin{eqn*}
\det
\Lambda_k^T=\prod_{j=2}^{\ell-k-1}\prod_{i=2}^{\frac{k}{2}}b_{\ell-j}b_{\ell-i}\det
\left(
\renewcommand{\arraystretch}{0.9}
\begin{array}{ccc}
b_{\ell-1}&b_{\ell}&0\\
b_{\ell-k-1}&0&b_{\ell}\\
0&b_{\ell-k-1}&b_{\ell-1}
\end{array}
\right).
\end{eqn*}%
Therefore  $\det
\Lambda_k^T=-2b_{\ell-1}b_{\ell}b_{\ell-k-1}\prod_{j=2}^{\ell-k-1}\prod_{i=2}^{\frac{k}{2}}b_{\ell-j}b_{\ell-i}$,
which is non-zero. We then  deduce  that the  rank of  $\Lambda_k^T$ is  $d_{k+1}$.

We  now consider the case where $k$ is odd. The root of $\g$ of  length
 $k$ are ordered in lexicographic order, to wit:
 $
\gamma_1=\l_1-\l_{k+1},
\dots,
\gamma_{\ell-k-1}=\l_{\ell-k-1}-\l_{\ell-1},
\gamma_{\ell-k}=\l_{\ell-k}-\l_{\ell},
\gamma_{\ell-k+1}=\l_{\ell-k}+\l_{\ell},
\gamma_{\ell-k+2}=\l_{\ell-k+1}+\ell_{\ell-1},
\dots,
\gamma_{\ell-\frac{k-1}{2}-2}=\l_{\ell-\frac{k-1}{2}-3}+\l_{\ell-\frac{k-1}{2}+2},
\gamma_{\ell-\frac{k-1}{2}-1}=\l_{\ell-\frac{k-1}{2}-2}+\l_{\ell-\frac{k-1}{2}+1},
\gamma_{\ell-\frac{k-1}{2}}=\l_{\ell-\frac{k-1}{2}-1}-\l_{\ell-\frac{k-1}{2}},
$
and the roots of $\g$ of length $k+1$ in lexicographic order, which gives the  array below where all the decompositions
of a root of length $k+1$ as a sum of a simple root with a root of length $k$ have been indicated on the right column:
\begin{eqn*}
\renewcommand{\arraystretch}{0.9}
\begin{array}{lllll}
\b_1&=&\l_1-\l_{k+2}&=&\left\{\begin{array}{l}
                                       \a_1+\gamma_2,\\
                                       \gamma_1+\a_{k+1},
                                        \end{array}
                                          \right.\\
\vdots&&\vdots&&\\
\b_{\ell-k-1}&=&\l_{\ell-k-1}-\l_{\ell}&=&\left\{\begin{array}{l}
                                       \a_{\ell-k-1}+\gamma_{\ell-k},\\
                                       \gamma_{\ell-k-1}+\a_{\ell-1},
                                        \end{array}
                                          \right.\\
\b_{\ell-k}&=&\l_{\ell-k-1}+\l_{\ell}&=&\left\{\begin{array}{l}
                                       \a_{\ell-k-1}+\gamma_{\ell-k+1},\\
                                       \gamma_{\ell-k-1}+\a_{\ell},
                                        \end{array}
                                          \right.\\
\b_{\ell-k+1}&=&\l_{\ell-k}+\l_{\ell-1}&=&\left\{\begin{array}{l}
                                       \a_{\ell}+\gamma_{\ell-k},\\
                                       \gamma_{\ell-k+1}+\a_{\ell-1},\\
                                        \gamma_{\ell-k+2}+\a_{\ell-k},
                                        \end{array}
                                          \right.\\
\b_{\ell-k+2}&=&\l_{\ell-k+1}+\l_{\ell-2}&=&\left\{\begin{array}{l}
                                       \a_{\ell-k+1}+\gamma_{\ell-k+3},\\
                                        \gamma_{\ell-k+2}+\a_{\ell-2},
                                        \end{array}
                                          \right.\\
\vdots      &&\vdots&& \\
\b_{\ell-\frac{k-1}{2}-2}&=&\l_{\ell-\frac{k-1}{2}-3}+\l_{\ell-\frac{k-1}{2}+1}&=&\left\{\begin{array}{l}
                                       \a_{\ell-\frac{k-1}{2}-3}+\gamma_{\ell-\frac{k-1}{2}-1},\\
                                        \gamma_{\ell-\frac{k-1}{2}-2}+\a_{\ell-\frac{k-1}{2}+1},
                                        \end{array}
                                          \right.\\
\b_{\ell-\frac{k-1}{2}-1}&=&\l_{\ell-\frac{k-1}{2}-2}+\l_{\ell-\frac{k-1}{2}}&=&\left\{\begin{array}{l}
                                       \a_{\ell-\frac{k-1}{2}-2}+\gamma_{\ell-\frac{k-1}{2}},\\
                                        \gamma_{\ell-\frac{k-1}{2}}+\a_{\ell-\frac{k-1}{2}-2}
                                        \end{array}
                                          \right.
\end{array}%
\end{eqn*}%
The matrix $^t\Lambda_k$, defined in (\ref{matrgamma}) has the following form:
\begin{eqn*}
^t\Lambda_k=\left(
\renewcommand{\arraystretch}{0.9}
\begin{array}{cccccccccc}
b_{k+1}&      &      & & & & & & & \\
b_1   &\ddots& & & & & & & & \\
      &\ddots&\ddots& & & & & & & \\
      &      &\ddots&b_{\ell-1}  &b_{\ell}& & & & & \\
      &      &      &b_{\ell-k-1}&0          &b_{\ell}& & & & \\
      &      &      &           &b_{\ell-k-1}&b_{\ell-1}& & & & \\
      &      &      &           &           &b_{\ell-k}&b_{\ell-2}& & &\\
      &      &      &           &           &         &b_{\ell-k+1}&\ddots&&\\
      &      &      &           &           &         &           &\ddots&b_{\ell-\frac{k-1}{2}+1}&\\
      &      &      &           &           &         &           &      &b_{\ell-\frac{k-1}{2}-3}&b_{\ell-\frac{k-1}{2}}\\
      &      &      &           &           &         &           &  &&b_{\ell-\frac{k-1}{2}-2}
\end{array}
\right).
\end{eqn*}%
By removing the first line of  $\Lambda_k^T$, we  obtain a upper  triangular square  ($d_{k+1}\times d_{k+1}$)
matrix   $\Gamma_k$ and which is of  rank $d_{k+1}$,  when
 $b_{k+1},\dots b_{\ell}$ are  all non-zero. This implies
that the rank of $\Lambda_k $ is $d_{k +1}$.\\
\emph{(b)} For the exceptional simple Lie algebras  $G_2,
  F_4,E_6,E_7$ and  $E_8$, we check the result by a direct computation on the software Maple.
We give  the program {\tt{Maple}} that completes the proof of Proposition  \ref{prkhim1}.
 We restrict ourself to the Lie algebra $E_6$ (for the other types, we use the same program with a adapted vector {\tt R}).\\
\smallskip
When  $\g$ is the simple Lie algebra of type $E_6$, the cardinal of the set of positive roots of $\g$ is  $N:=36$. 
We suppose that the elements of  $\Phi_+$ are indexed by lexicographic order. 
To each  $\a$ of  $\Phi_+$, we associate a row vector $R[i]:=[a_1,\dots,a_6]$ such that $\a=\sum_{j=1}^{6}a_j\a_j$,
where $\a_1, \dots, a_\l$ are the simple roots.\\
\smallskip
\noindent
{\tt with(linalg):}\\
 \noindent
{\tt N:=36:}\\
\noindent
{\tt rank:=6;}\\
\noindent
\begin{tabular}{lll}
\tt R[1]:=\tt [1,0,0,0,0,0]:&
\tt R[2]:=\tt [0,1,0,0,0,0]:&
\tt R[3]:=\tt [0,0,1,0,0,0]:\\
\tt R[4]:=\tt [0,0,0,1,0,0]:&
\tt R[5]:=\tt[0,0,0,0,1,0]:&
\tt R[6]:=\tt[0,0,0,0,0,1]:\\
\tt R[7]:=\tt[1,0,1,0,0,0]:&
\tt R[8]:=\tt[0,1,0,1,0,0]:&
\tt R[9]:=\tt [0,0,1,1,0,0]:\\
\tt R[10]:=\tt[0,0,0,1,1,0]:&
\tt R[11]:=\tt[0,0,0,0,1,1]:&
\tt R[12]:=\tt[1,0,1,1,0,0]:\\
\tt R[13]:=\tt[0,1,1,1,0,0]:&
\tt R[14]:=\tt[0,1,0,1,1,0]:&
\tt R[15]:=\tt [0,0,1,1,1,0]:\\
\tt R[16]:=\tt [0,0,0,1,1,1]:&
\tt R[17]:=\tt[1,1,1,1,0,0]:&
\tt R[18]:=\tt[1,0,1,1,1,0]:\\
\tt R[19]:=\tt[0,1,1,1,1,0]:&
\tt R[20]:=\tt[0,1,0,1,1,1]:&
\tt R[21]:=\tt[0,0,1,1,1,1]:\\
\tt R[22]:=\tt[1,1,1,1,1,0]:&
\tt R[23]:=\tt[0,1,1,2,1,0]:&
\tt R[24]:=\tt[1,0,1,1,1,1]:\\
\tt R[25]:=\tt[0,1,1,1,1,1]:&
\tt R[26]:=\tt[1,1,1,2,1,0]:&
\tt R[27]:=\tt[1,1,1,1,1,1]:\\
\tt R[28]:=\tt[0,1,1,2,1,1]:&
\tt R[29]:=\tt[1,1,2,2,1,0]:&
\tt R[30]:=\tt[1,1,1,2,1,1]:\\
\tt R[31]:=\tt[0,1,1,2,2,1]:&
\tt R[32]:=\tt[1,1,2,2,1,1]:&
\tt R[33]:=\tt[1,1,1,2,2,1]:\\
\tt R[34]:=\tt[1,1,2,2,2,1]:&
\tt R[35]:=\tt[1,1,2,3,2,1]:&
\tt R[36]:=\tt[1,2,2,3,2,1]:
\end{tabular}%

\noindent
\#{\tt We define a procedure to calculate the length of 
  a
root X}\\
\newcommand{\es}{\hspace*{0.5cm}}
\noindent
{\tt long:=proc(X)}  \\
{\tt \es sum(X[k],k=1..nops(X))}\\
 {\tt $\quad$end:}\\
\noindent
\#{\tt We construct a list containing the roots of the same length}\\
\noindent
{\tt lis:=}{\tt proc(i)}\\
{\tt local k, list;}\\
{\tt  list:=\tt [];}\\
{\tt  \es for k from 1 to N do}\\
{\tt   \es\es if long(R[k])=i then} \\
{\tt \es\es\es  list:=[op(list),R[k]]}\\
  {\tt \es \es fi}\\
 {\tt \es od}\\
{\tt end:}\\
\noindent
\#{\tt  Relation between  a root  i  of length
k  and   a root  j  of length   k+1}\\
\noindent
{\tt a:=\tt proc(k,i,j)}\\
\noindent {\tt local l,res,dL;}\\
\noindent {\tt res:=0:}\\
\noindent {\tt dL:=lis(k+1)[j]-lis(k)[i];}\\
\es {\tt for l from 1 to rank do}\\
 \es\es {\tt if dL=R[l] then res:=b[l]}\\
 \es \es {\tt fi;}\\
\es {\tt od;}\\
{\tt res}\\
{\tt end:}\\
\noindent
{\tt Gammas:=proc(k)}\\
 {\tt\es  matrix(nops(lis(k)),nops(lis(k+1)),(i,j)->a(k,i,j))}\\
{\tt end:}\\
\noindent
\# {\tt We verify if the rank of the matrix  $\Gamma_k$ (that is
$\Lambda_k^T $ in the proof  of \\
\noindent
 \# Proposition \ref{prkhim1})  is the number of roots of length $k$.}\\
\noindent
{\tt verif:=proc(k)}\\
 {\tt \es if nops(lis(k+1))-rank(Gammas(k))=0 then  1 else 0}\\
 {\tt \es fi}\\
{\tt end:}\\
\noindent
{\tt K:=1;}\\
{\tt for i from 1 to long(R[N])-1 do  K:=K*verif(i)} \\
{\tt od:}\\
 \noindent
{\tt if K=1  then print(OK)} \\
  {\tt \es else print("pas OK")} \\
{\tt fi;}\\
\es\es\es\es \es\es\es\es\es\es {\tt OK}
\end{proof}
\section{A conjectured integrable system}\label{sect4}
We believe that the periodic Full Kostand-Toda lattice and the periodic Toda lattice are two extremes cases of a
 string of integrable systems, that we now present.
In Proposition \ref{fgyu}, we have shown that $\TT_\l$ is a Poisson submanifold of $\tilde{\g}$,
using the fact, stated in (\ref{rel_ps}), that
$$\TT_{\lambda}:=\bigoplus_{-m_\ell\leqslant i\leqslant 0}{\tilde{\g_i}}+f,$$
where $f:=\sum_{i=1}^{\ell}e_i+\lambda e_{-\beta}\in{\tilde \g}_1$.
The same argument shows that $ \TT_\l^{(k)}:= \bigoplus_{0 \leqslant i\leqslant k}{\tilde{\g}_{-i}}+f$
is a Poisson submanifold of $\tilde{\g}$ for all $k=1, \dots,m_\ell$.

By construction, the phase spaces $\TT_\l^{(m_{\ell})} $ and $\TT_\l^{(1)}$ are the phase spaces of the
periodic Full Kostand-Toda lattice and the periodic Toda lattice respectively.
Since the differential equation associated to the Hamiltonian $\frac12{\inn{x(\l)}{x(\l)}_\l} $
 is Liouville integrable  in the two extreme cases, it is natural to ask whether it is Liouville integrable for all $k$.
\\ More precisely, it is natural to ask whether the following differential equation is Liouville integrable for all $k=1,...,m_\ell$: 
\begin{equation}
\dot{L}^{(k)}(\l)=[L^{(k)}(\l),L^{(k)}(\l)_-], \forall 1<k<m_{\ell},
\end{equation}%
where $L^{(k)}(\l)$ is an element of the phase space 
\begin{equation}\label{phase-space-Tk}
\TT_{\l}^{(k)}:=\l e_{-\b}+\h+\sum_{1\leqslant j\leqslant k}\g_{-j}+\l^{-1}\g_{m_{\ell}+1-j}
\end{equation}%
and ${L^{(k)}(\l)}_-$ is the strictly lower part of $L^{(k)}(\l)$.
\begin{example}
When $\g$ is $\Liesl{n}(\C)$ and $\h$ is the subalgebra of diagonal matrices of $\Liesl{n}(\C)$,   
an element $L^{(k)}(\l)$ 
of $\TT_{\l}^{(k)}$ has the following form:
\begin{equation}
\renewcommand{\arraystretch}{0.3}
\left(
\begin{array}{ccccccccc}
a_{11}   &1     &0     &\ldots&0      &\l^{-1}a_{1,n-k+1}&\ldots &\ldots      &\l^{-1}a_{1,n}\\
\vdots  &\ddots&\ddots&\ddots&       &\ddots          &\ddots &            &\vdots\\
\vdots  &      &\ddots&\ddots&\ddots &                &\ddots &             &\l^{-1}a_{k,n}\\
\vdots  &      &      &\ddots&\ddots &\ddots          &      &             &0\\
a_{k+1,1}&      &      &      &\ddots &\ddots          &\ddots&             &\vdots\\
0       &\ddots&      &      &       &\ddots          &\ddots& \ddots       &\vdots\\
\vdots  &\ddots&\ddots&      &       &                &\ddots&\ddots        &0\\
0       &      &\ddots&\ddots&       &                &      &\ddots        &1\\
\l      &0     &\ldots&0     &a_{n,n-k}&\ldots          &\ldots&\ldots        &a_{nn}
\end{array}
\right).
\end{equation}
Notice that these differential equations are those that appear in \cite{MMD}, for formal solutions are given.
\end{example}
For the family of functions that give the Liouville integrability, there is again a natural candidate,
given by the restriction of the family $(\tilde{F}_{i,j}, 1\leqslant i\leqslant \ell, 0\leqslant i\leqslant m_{\ell})$ to
 $\TT_{\l}^{(k)}$. Again, several of these restrictions
 vanish or are constant. It seems that the following families of functions:
  $$ {\tilde F}^{(k)}=(\tilde{F}_{j,i}, ~~1\leqslant i\leqslant, ~~ 1\leqslant j\leqslant  E( k \frac{m_i+1}{m_\ell+1}) $$
admit a restrictions to $\TT_\l^{(k)} $ which are independent. At least, we have been able to check,
 with Maple, that these restrictions are independent for $\Liesl{n}(\C)$ with $n=2, \dots,7$, 
and for the Lie algebras $B_n$ for $n=2, \dots,6$ for all possible value of $k$.
 For all the previous cases, we have also verified, by using Maple, that
 the rank $\Rk(\TT_{\l}^{(k)}, {\PB}_{\tilde R})=\dim\TT_{\l}^{(k)}-\frac12\Rk(\TT_{\l}^{(k)},\PB_{\tilde R})$
 of the restricted Poisson structure satisfies the third item of Definition \ref{def:liouville},
 which establishes the Liouville integrability.
We therefore think that this should be always true.
\begin{conjecture}
  The triplet $(\TT^{(k)}_{\l},{\tilde\FF}^{(k)}_{\lambda},{\PB}_{\tilde R})$ is an integrable system.
\end{conjecture}
The first difficulty is that, for $1< k < m_{\ell}$, it is not possible any more to find 
in the phase space of $\TT_\l^{(k)}$ points where we can apply  Theorem  \ref{indvki} of Ra\"is:
 we therefore probably have to find a suitable generalization of this result. It is very likely
 that we have to use a point of the form  
$$L_0(\l)=\l e_{-\b}+e +\sum_{i=1}^{\ell}b_ih_i +\sum_{i=1}^{d_k}a_i e_{-\gamma_i}+\l^{-1}\sum_{i=1}^{d_{m_{\ell}+1-k}}c_i e_{\eta_i},$$
 where 
   $\gamma_1,\dots, \gamma_{d_k}$ are the $d_k$ roots of $\g$ of   length $k$ and $\eta_1, \dots,\eta_{d_{m_{\ell}+1-k}}$ are the 
$ d_{m_{\ell}+1-k}$ roots of $\g$ of length $m_{\ell}+1-k$.
 Also, it is not clear to see at which point one should compute the rank. It is even far from being easy to 
guess which ones of the functions ${\tilde\FF}^{(k)}_{\lambda}$ are going to be Casimir functions.
 It is clear that only the functions $ \tilde{F}_{E( k \frac{m_i+1}{m_\ell+1} ),i}$ 
may be Casimirs functions, but some of them are not. For instance, for $k=1$, only one of them
 (for $i=\ell, j= k$) is a Casimir function, while for the periodic Full Kostant-Toda, all the functions 
$ F_{m_i,i}$ for $i=1, \dots,\ell$ are Casimirs (by Proposition \ref{caspft}). For generic values of $k$, 
the behavior seems at first to be quite random. For instance, in the case $\g= \Liesl{n}(\C)$, $n=7$  
and $k=2$, respectively $k=3$, (cases where the Liouville integrability can be proved by Maple), the Casimir functions
 are $ F_{3,1}, F_{6,2}  $, respectively $F_{2,1}, F_{4,2},F_{6,3} $.      

%


%
%

\end{document}